\documentclass[12pt]{amsart}
\usepackage{amssymb,amscd}
\usepackage[dvipdfmx]{graphicx,color}

\usepackage{graphicx}

\newtheorem*{Whitney towers}{Theorem~\ref{Whitney towers}}
\newtheorem*{h-towers}{Theorems ~\ref{half} \& \ref{$(n)$-solvable}}

\newtheorem*{surgery curves}{Theorem~\ref{surgery curves}}
\newtheorem*{cg=0}{Theorem~\ref{vanish}}

\newtheorem{thm}{Theorem}[section]
\newtheorem{pr}[thm]{Proposition} 
\newtheorem{lem}[thm]{Lemma}
\newtheorem{cor}[thm]{Corollary}

\newtheorem{cla}[thm]{Claim}

\theoremstyle{definition}
\newtheorem{defn}[thm]{Definition}
\newtheorem{que}[thm]{Question}

\newtheorem{note}[thm]{Note}

\numberwithin{equation}{section}
\numberwithin{figure}{section}

\newcommand{\x}{\times}

\newcommand{\Z}{\mathbb{Z}}
\newcommand{\N}{\mathbb{N}}
\newcommand{\C}{\mathbb{C}}
\newcommand{\Q}{\mathbb{Q}}

\newcommand{\R}{\mathbb{R}}

\def\yen{{\setbox0=\hbox{Y}Y\kern-.97\wd0\vbox{hrule height.lex width.98%
\wd0\kern.33ex\hrule height.lex width.98\wd0\kern.45ex}}}

\begin{document}
\title{
Brieskorn submanifolds, Local moves on knots, and knot products
}
\author{Louis H. Kauffman and  Eiji Ogasa}


\begin{abstract}
We first prove the following: 
Let $p\geq2$  
and $p\in\N$. 
Let $K$ and $J$ be closed, oriented, $(2p+1)$-dimensional 
connected, $(p-1)$-connected, simple submanifolds of $S^{2p+3}$. 
Then $K$ and $J$ are isotopic 
if and only if 
a Seifert matrix associated with a simple Seifert hypersurface for $K$ 
is $(-1)^p$-$S$-equivalent to 
that for $J$. 
We also discuss the $p=1$ 
case. 

This result implies one of our main results: 
 Let $\mu\in\N$. A 1-link $A$ is pass-equivalent to a 1-link $B$ if and only if 
$A\otimes^{\mu} {\rm Hopf}$ is 
$(2\mu+1, 2\mu+1)$-pass-equivalent to 
$B\otimes^{\mu} {\rm Hopf}$.

It also implies the other of our main results: 
We strengthen the authors' old result that   
two-fold cyclic suspension commutes with the performance of the twist move
for spherical $(2k+1)$-knots. See the body for the precise statement.  

Furthermore it implies the following:  
Let $p\geq2$ 
and $p\in\N$. Let $K$ be a closed oriented $(2p+1)$-submanifold of $S^{2p+3}$. Then 
$K$ is a Brieskorn submanifold if and only if 
$K$ is connected, $(p-1)$-connected, simple and 
has a $(p+1)$-Seifert matrix associated with a simple Seifert hypersurface 
 that is $(-1)^p$-$S$-equivalent to 
a $KN$-type 
(see the body of the paper for a definition). 
We also discuss the $p=1$  
case.  
\end{abstract}

\thanks{\hskip-4mm E-mail: kauffman@uic.edu\quad  
ogasa@mail1.meijigakuin.ac.jp \newline
\newline MSC2000 57Q45, 57M25, 32S55}

\date{}
\maketitle   

\tableofcontents

\bigbreak
\section{Introduction}\label{Introduction}
\noindent
We begin this introduction with the two basic theorems 
(Theorems \ref{Chicago} and \ref{carrot})
that we prove in this paper. 
Theorem \ref{Chicago} tells us how the Seifert matrix for a simple codimension two 
submanifold of an odd-dimensional sphere determines the embedding of the submanifold. 
 We then explain briefly, in this introduction, our applications of these results to knot products and local moves. 
The applications are the main results of this paper.   (Theorems \ref{Bunkyo} and \ref{aoiro}). 
\\
 
 The reader should recall that a codimension two submanifold $K$ of an odd-dimensional sphere is said to be {\it simple} if the fundamental group of the complement is $\Z$ and the higher homotopy groups of the complement vanish below the middle dimension of $K$. 
The Seifert pairing is defined on pairs of cycles on a manifold  $F$ whose boundary is $K.$ 
The cycles must have dimensions that add to the dimension of $F$, and the pairing is defined by taking the linking number of one cycle, with the other one pushed off the manifold $F$ into the positive normal direction. 
\\

A Brieskorn manifold $\Sigma(a_0,a_1,\cdots , a_n)$ is the intersection of 
 the complex variety of $z_{0}^{a_0} + \cdots + z_{n}^{a_{n}}$ with the unit sphere in complex $n+1$-space. 
Here, we obtain a smooth submanifold in the unit sphere. We call it a Brieskorn submanifold.   
Brieskorn (sub)manifolds have been a rich source of examples for knots in higher dimensions and exotic differentiable structures that are related to the knot theory of a codimension-two embedding. In his seminal book \cite{Milnor} ``Singular points of complex hypersurfaces" John Milnor put the subject of the Brieskorn manifolds in the context of links of complex hypersurface singularities and proved a fundamental fibration theorem. 
\\

Many structures are related across dimensions in these families of algebraic varieties. We say that a knot is in dimension $n$ if it is an $n-2$ dimensional submanifold of the sphere $S^{n}$. Thus a classical knot is in dimension $3$.
In earlier papers \cite{Kauffman, KauffmanNeumann} we  have studied constructions such as the {\it knot product} that generalizes the link of the sum of two singularities to a construction that produces a new knot in dimension $n+m+1$ from knots in dimensions $n$ and $m$ 
 respectively (well defined when one knot is fibered). This generalization is significant because it utilizes structure that formerly was only available through algebraic varieties to knot theory proper. 
\\

In more recent work  (\cite{KauffmanOgasa, KauffmanOgasaII}), 
the authors of the present paper have worked with other constructions  
(local moves on $n$-knots for any natural number $n$, and generalized band-passing)    
in relation to the knot product construction.  
In all these cases it has been important to study the Seifert pairing for manifolds bounding the knots that we study. 
One of our results there: 
If a 1-knot $J$ is obtained from another 1-knot $K$,  
then the knot product of $J$ and the Hopf link, which is a 3-knot, 
is obtained from that of $K$ and it by one twist move, 
which is a local move on 5-knot. 
\\

In this paper we prove theorems showing how, for simple submanifolds, the Seifert pairing can classify the submanifold, and how under certain conditions, knowing the Seifert pairing characterizes the manifold as a Brieskorn manifold. We believe that these theorems are satisfying, and that they lend better understanding to this subject of codimension two submanifolds of spheres.

\def\iruka
{Let $A$ be an $r\x r$-matrix $(r\in\N\cup\{0\}$).    
(Let $\N$ be the set of natural numbers. Note that $n$ is a natural number if and only if $n$ is a positive integer.)   
We give two basic forms of {\it $(-1)^p$-$S$-equivalence}.    
We say $A \sim A'$ 
if \newline
 $A'=$
$\begin{pmatrix}
A&\alpha&0\\
(-1)^p\hskip1mm^t\hskip-1mm\alpha&0&0\\
0&1&0
\end{pmatrix}$ 
(respectively,  
$\begin{pmatrix}
A&\beta&0\\
(-1)^p\hskip1mm^t\hskip-1mm\beta&0&1\\
0&0&0
\end{pmatrix}$), 
where $\alpha$ and $\beta$ are column vectors, 
or if $A'$ is equivalent to $A$ 
(that is, if $A'$ is a conjugate of $A$, or if $A'=^t\hskip-2mm PAP$).   
Then we say that 
$(-1)^p$-$S$-equivalence is the equivalence relation generated by these basic equivalences.}

\def\odoroita{
$(1)$ 
Let $p\geq2$  
and $p\in\N$. 
Let $K$ and $J$ be closed oriented \newline
$(2p+1)$-dimensional connected, $(p-1)$-connected, simple submanifolds of $S^{2p+3}$. 
Then the following two statements are equivalent.     

\smallbreak\noindent{\rm(i)}   
$K$ is isotopic to $J$. 

\smallbreak\noindent{\rm(ii)}  
A simple 
Seifert matrix $P_K$ for $K$ is $(-1)^p$-$S$-equivalent to 
a simple 
Seifert matrix $P_J$ for $J$.

\bigbreak\noindent$(2)$ 
Let $K$ and $J$ be closed oriented 3-dimensional simple submanifolds of $S^5$. 
Then the following two statements are equivalent.     

\smallbreak\noindent{\rm(i)}   
$K$ is isotopic to $J$.

\smallbreak\noindent{\rm(ii)}  
There is a simple Seifert hypersurface $V_K$ $($respectively,  $V_J)$ for $K$ $($respectively,  $J)$ 
with the following properties: 
There is an orientation preserving diffeomorphism map $f:V_K\to V_J$. 
  $V_K$ $($respectively, $V_J)$ consists of one 0-handle $h^0_{K}$ $($respectively,  $h^0_{J})$ 
and 2-handles $h^2_{Ki}$ $($respectively,  $h^2_{Jj})$, 
where $i,j\in\{1,...,\nu\}$ and $\nu$ is a nonnegative integer 
$($Note that if $\nu=0$, we regard the set as the empty set$)$.  
$f(h^0_{K})=h^0_{J}$. 
$f(h^2_{Kl})=h^2_{Jl}$ for each $l\in\{1,...,\nu\}$. 
$s(f(h^2_{Ka}), f(h^2_{Kb}))=s(h^2_{Ja}, h^2_{Jb})$ for each pair 
$(a,b)\in\{1,...,\nu\}\x\{1,...,\nu\}$, 
where $s(h^2_{*a}, h^2_{*b}) (*=K, J)$ denotes a Seifert pairing of 
a pair of 2-cycles which are defined by $h^2_{*\natural} (\natural=a,b)$. 
}


\bigbreak
The following theorems 
(Theorems  \ref{Chicago} and \ref{carrot}) 
are fundamental to the structure of simple submanifolds of spheres and have not been proven before. So we prove them in this paper.  
We apply these results to obtain proofs of our main new results, 
Theorems \ref{Bunkyo} and \ref{aoiro},  
about knot products and local moves on high dimensional knots.

Mainly in \S\ref{prod}, \ref{move}, and \ref{RkB}  are 
the other terms and definitions than those which are explained in this section  but 
which are needed for these theorems.

\begin{defn}\label{saisho}
\iruka 
\end{defn}

\begin{thm}\label{Chicago}
\odoroita 
\end{thm}

\def\honmaya
{$(1)$ 
Let $p\geq2$  
and $p\in\N$. Let $K$ be a closed oriented $(2p+1)$-submanifold of $S^{2p+3}$. 
Then the following two statements are equivalent.     

\smallbreak\noindent{\rm(i)} 
$K$ is a Brieskorn submanifold. 

\smallbreak\noindent{\rm(ii)}  
$K$ is connected, $(p-1)$-connected, simple, and 
has a simple Seifert matrix 
$P$ that is $(-1)^p$-$S$-equivalent to a $KN$-type.  
$($See Definition \ref{KN} for the definition of  $KN$-types.$)$

\bigbreak\noindent$(2)$ 
 Theorem $\ref{Chicago}.(2)$ is true  
if $J$ is the Brieskorn submanifold which is a 3-manifold. 
%
%
%
%
}


Theorem \ref{Chicago} and Definition \ref{KN} imply Theorem \ref{carrot}.

\begin{thm}\label{carrot}  
\honmaya 
\end{thm}


These theorems allow us to make comparisons of the action of 
the knot product 
(such as taking a knot product with a Hopf link) 
and certain local moves on high dimensional knots. 
In Theorem \ref{Bunkyo}, 
we strengthen the authors' old results in \cite{KauffmanOgasa, KauffmanOgasaII}  
that 
taking a knot product with the Hopf link commutes with the performance of the pass move.   
In Theorem \ref{aoiro}, 
we strengthen the authors' old results  in \cite{KauffmanOgasa, KauffmanOgasaII} 
that 
two-fold cyclic suspension commutes with the performance of the twist move. 
Much of the work in this paper depends upon familiarity with knot products and 
local moves on knots. 
We take care to review these matters in the next three sections.



\bigbreak
\section{Review of  knot products}\label{prod}
\noindent 
In 
\cite{KauffmanOgasa} we began to research relations between knot products and 
local-moves on knots. 
This paper is a sequel to 
\cite{KauffmanOgasa, KauffmanOgasaII}. 
Knot products, or products of knots, were defined and have been researched in 
\cite{Kauffman, KauffmanNeumann, KauffmanOgasa}. 
Local moves on high dimensional knots were defined and 
have been  researched 
in 
\cite{KauffmanOgasa, Ogasa98n,  Ogasa02,  Ogasa04, Ogasa07,   Ogasa09, OgasaT3, OgasaIH}. 
\cite{KauffmanOgasaB} is a preprint of this paper. 
\bigbreak

We work in the smooth category unless we indicate otherwise.

Let $n, m\in\N$.  
An  {\it $($oriented$) ($ordered$)$ $m$-component $n$-$($dimensional$)$ link} 
 is a smooth, oriented submanifold $L=(L_1,...,L_m)$ $\subset$ $S^{n+2}$, 
which is the ordered disjoint union of $m$ manifolds, each PL homeomorphic 
to the standard $n$-sphere. If $m=1$, then $L$ is called a {\it knot}.

\bigbreak 
Note the following: 
We usually define $n$-links as above (see e.g. 
\cite{CochranOrr}). 
Not all $n$-knots are diffeomorphic to the standard $n$-sphere 
although all $n$-knots are PL homeomorphic to the standard $n$-sphere.  
The reason for this is the fact that many exotic $n$-spheres,  
which are not diffeomorphic to the standard $n$-sphere,  
can be embedded smoothly in $S^{n+2}$ 
(see 
\cite{Levinecob, Levinesimp, Milnor} for the proof of this fact.) 
\bigbreak

%
%
%
%
%
%
%
%
Let $n$ be a positive integer. 
Two submanifolds $J$ and  $K$ $\subset S^n$ are {\it $($ambient$)$ isotopic}  
if there is a smooth orientation preserving family of diffeomorphisms $\eta_t$ of $S^n$, $0\leqq t\leqq1$, with $\eta_0$ the identity and $\eta_1(J)=K$.   
An $m$-component $n$-link $L=(L_1,...,L_m)$ is called the {\it trivial $(n$-$)$link}   
if each $L_i$ bounds an $(n+1)$-ball $B_i$ 
embedded in $S^{n+2}$ 
and if $B_i\cap B_j=\phi$  
for each distinct $i, j$. 
If $L$ is the trivial 1-component ($n$-)link,  then $L$ is called the {\it trivial $(n$-$)$knot}. 

\bigbreak 
Note the following: 
When we study singular points of complex hypersurfaces 
(see e.g. 
\cite{Milnor}) and 
knot products (see e.g. 
\cite{Kauffman, KauffmanNeumann}),   
`knots' often mean (not necessarily connected or PL spherical) 
$n$-dimensional closed submanifolds of $S^{n+2}$.  
In these cases we use the following terms for clarity. 
We say that an $n$-submanifold $M$ is 
a {\it PL spherical knot}  (respectively,  {\it topological spherical knot})  
if $M$ is PL homeomorphic (respectively,  homeomorphic) to the standard $n$-sphere.  
Recall the following well-known fact: Let $n\neq4$. 
$M$ is homeomorphic to  the standard $n$-sphere if and only if 
$M$ is PL homeomorphic to  the standard $n$-sphere. 
If $n=4$, it is an open problem whether this fact holds. 

So, in the $n\neq4$ case,  
we say that $K$ is a {\it spherical knot} 
if $K$ is a PL spherical knot (respectively,  topological spherical knot). 
Of course we do not use these terms when it is clear from the context. 

\bigbreak  
Let $l\in\N.$ 
Let $L$ be a (not necessarily connected or PL spherical) closed oriented 
$l$-submanifold of $S^{l+2}$. 
Let $W$ be a connected compact oriented $(l+1)$-submanifold of $S^{l+2}$ 
such that $\partial W=A$. 
We call $W$ a {\it Seifert hypersurface} for $A$.   
By obstruction theory, such $W$ exists. See e.g. 
\cite[Theorem 3 in P.50]{Kirby}. 
Replace $S^{n+2}$ with another (not necessarily closed) compact oriented $(n+2)$ manifold. 
If there is such $W$ for $L$, we also call $W$ a {\it Seifert hypersurface} for $L$.

In this paper we abbreviate 
manifold-with-boundary  (respectively,  submanifold-with-boundary)  
to 
manifold (respectively,  submanifold) when this is clear from the context.

\bigbreak

We review the definition of knot products, or products of knots, after we explain a question 
which motivates the definition of knot products.  
This question was discussed and answered in \cite[\S 4]{KauffmanNeumann}. 

Let $a$ be a positive integer. 
Let $f: \C^a \longrightarrow\C$ 
be a (complex) polynomial mapping with an isolated singularity at the origin of $\C^a$
and satisfy $f$(the origin)=the origin. 
Recall that  one defines 
the {\it link of this singularity}, $L(f)=f^{-1}(0)\cap S^{2a-1}\subset S^{2a-1}.$ 
(See  \cite{Milnor} for the precise definition of the link of singularity.)
Here the symbol $f^{-1}(0)$ denotes the variety of $f$,  
and $S^{2a-1}$ denotes a sufficiently small sphere about the origin of $\C^a.$

Let $b$ be a positive integer. Given another such polynomial $g:\C^{b} \longrightarrow \C$, 
 form $f + g: \C^{a+b} =\C^a \times \C^b\to \C$ 
by $(f+g)(x,y)=f(x)+g(y)$. 

What kind of relation do we have among  
$L(f) \subset S^{2a-1}$, $L(g) \subset S^{2b-1}$, and $L(f + g)\\\subset S^{2a+2b-1}$?  
An answer to this question is explained 
below in the form of a construction that works for more general codimension two embeddings. 
This construction is due to \cite{KauffmanNeumann}. 
\\

Knot products are defined in \cite[Definition 3.3]{KauffmanNeumann}. 
An alternative definition is given in \cite[Lemma 3.4]{KauffmanNeumann}, 
which we cite almost verbatim below. 
We emphasize that the reference  \cite{KauffmanNeumann} contains the proofs about 
the well-definedness of this construction.

A knot $(S^k, K)$ 
is defined in \cite[the fourth paragraph of page 371]{KauffmanNeumann} as   
a (not necessarily connected or spherical) $(k-2)$-dimensional closed oriented submanifold $K$ of $S^k$.

A fibered knot $(S^l, L, b)$ 
is defined in \cite[the fourth paragraph of page 371]{KauffmanNeumann} as   
a (not necessarily connected or spherical) $(l-2)$-dimensional fibered closed oriented submanifold $L$ of $S^l$.

The key to the knot product construction is the association of maps to $D^2$ (as explained below) 
with a fibered closed oriented submanifold $(S^l,L,b)$. 
The fibered submanifold has a map $b:\overline{S^l-N(L)}\to S^1$
that is a fibration and restricts to the projection to $S^1$ on the boundary 
of $N(L)\cong L\x D^2$, where $\cong$ denotes an orientation preserving diffeomorphism. 
Thus $b$ extends to a map (still called) $b:S^l\to D^2$ with $b^{-1}(0)=L$. We can then 
form the cone on this map, $cb:D^{l+1}\to D^2$ with $(cb)^{-1}(0)=CL$, 
where $CL$ denotes the cone on $L$. 
Define a pair $(D^{2a}, f^{-1}(0))$ by using the algebraic variety several paragraphs above, 
where $D^{2a}$ is a $2a$-ball differentiably embedded in $\C^{2a}$ whose boundary is $S^{2a-1}$. 
The pair $(D^{l+1}, CL)$ is 
a generalization of $(D^{2a}, f^{-1}(0))$   
to the  case where we do not use an algebraic variety. 
\\

\noindent{\bf Definition of knot products.} 
Let $k$ and $l$ be positive integers.  Let $S^k=\partial D^{k+1}$ and $S^l=\partial D^{l+1}$.

\noindent{\bf Lemma 3.4 of \cite{KauffmanNeumann}.} 
{\it 
Let $(S^k, K)$ be a knot 
and 
$(S^l,L,b)$ a fibered knot. 
Let $F\subset D^{k+1}$ be a spanning manifold for $K$ 
as 
in \cite[Definition 3.1]{KauffmanNeumann}. 
Use \cite[Lemma 2.3]{KauffmanNeumann} to obtain $\gamma:D^{k+1}\to D^2$ with 
$\gamma^{-1}(0)=F,0$ a regular value of $\gamma$. 
Let $\tau:D^{l+1}\to D^2$ be a smoothing of $cb$.  
Use these maps to form the pullback 

$$
\begin{CD}
b(D^{k+1},F)&@>>>&D^{l+1}\\
@VVV& &@VV{\tau}V\\
D^{k+1}&@>{\gamma}>>&D^2
  \end{CD}
$$
\\


\noindent 
Thus 
$b(D^{k+1},F)\subset D^{k+1}\x D^{l+1}$. 
Then the knot product $(S^{k+l+1}, K\otimes b)$ 
is obtained by taking boundaries from this embedding. 
That is,   
$(S^{k+l+1}, K\otimes b)$ is isotopic to \\
$(\partial(D^{k+1}\x D^{l+1}), \partial(b(D^{k+1},F)))$.}

$(S^{k+l+1}, K\otimes b)$  is well-defined in the differentiable category in terms of the embedding $K$ and the fiber structure of $L$. 
\\

For 
$(S^k, K)$ and 
$(S^l, L, b)$, 
$(S^{k+l+1}, K\otimes b)$ is indicated by $(S^{k+l+1}, K\otimes L)$ as written in  
\cite[the first paragraph of \S 3]{KauffmanNeumann}. 



See \cite[Fig. 1 in page 371 and the part above it]{KauffmanNeumann} for 
the word `smoothing' in the cited lemma. 
In particular, a smoothing of $cb:D^{l+1}\to D^2$ (as in \cite[Lemma 3.4]{KauffmanNeumann}) 
with $(cb)^{-1}(0)=CL$ is obtained by changing to a nearby smooth map 
$\widetilde{cb}:D^{l+1}\to D^2$ with $(\widetilde{cb})^{-1}(0)=G$, 
where $G\subset D^{l+1}$ is a smooth submanifold of $D^{l+1}$ 
and $(\widetilde{cb})^{-1}(0)\cap S^{l}$ is isotopic to $L$. 
The manifold $G$ is part of a family of submanifolds 
that deforms to the singular subspace $CL$.

In this paper we let $K\otimes L$ denote 
$(S^{k+l+1}, K\otimes L)$ in order to make the notation shorter.  
Note that in 
\cite[the first paragraph of \S 3, and Definition 3.1]{KauffmanNeumann},  
$K\otimes L$ denotes the diffeomorphism type of the submanifold $(S^{k+l+1}, K\otimes L)$, 
$K\otimes L\hookrightarrow S^{k+l+1}$.   

Note that there are infinitely many cases 
a knot product $K\otimes L$ is defined 
even if neither $K$ nor $L$ is a spherical knot. (Recall that any spherical knot is connected.)

If $k$ (respectively, $l$) is one, we regard $K$ (respectively $L$) as the empty knot, 
which is defined in \cite[line($-3$) of page 371]{KauffmanNeumann}. 
The empty knots $[a]$ correspond to the maps \\$\tau_a:D^2\to D^2, \tau_a(z)=z^a$, 
where $a$ is an integer$\geqq2$. 
\cite[Lemmas 5.3 and 5.4]{KauffmanNeumann} 
explains taking a knot product with the empty knot [2]. 
It is an example of knot products.

$b(D^{k+1},F)$ is a submanifold of $D^{k+1}\x D^{l+1}$ 
as \cite[the first line of page 374]{KauffmanNeumann} is pointed out:  
In \cite[the last line and the diagram on it of page 373]{KauffmanNeumann},   
let $M$ be $D^{k+1}$, $D^{n+1}$ be $D^{l+1}$ and $V$ be $F$  
(note this $F$ is not `$F$ in \cite[Definition 2.1]{KauffmanNeumann}').   
Then $N=D^{k+1}\x D^{l+1}$. 
%
Note  
$N=\tau(M, \alpha)$ in \cite[the third line of Theorem 2.2]{KauffmanNeumann}, 
$\tau(M, \alpha)=b(M, \alpha)$ in \cite[Theorem 2.2.(ii)]{KauffmanNeumann}, and    
$b(M,\alpha)=b(M,V)$ in  \cite[Remark in page 376]{KauffmanNeumann}.  
Since $b(D^{k+1},F)$ is an inverse image of a regular value by a smooth map,  
$b(D^{k+1},F)$ is a submanifold of $D^{k+1}\x D^{l+1}$.

Hence $\partial(b(D^{k+1},F))$ 
is a well-defined differentiable submanifold of $S^{k+l+1}$.

The first remark in \cite[page 380]{KauffmanNeumann} explains the reason for the following: 
When we consider the knot product of two submanifolds,  
we impose the condition that one of the two submanifolds is fibered. 
If neither submanifold is fibered, the construction is not uniquely defined. 
\\

\cite[Corollary 3.6]{KauffmanNeumann} implies that for fibered knots $(S^l,L,b)$ and $(S^{l'},L',b')$,  \\
$(S^{l+l'+1}, L\otimes L')$ is isotopic to $(-1)^{(l-1)(l'-1)}(S^{l+l'+1}, L'\otimes L)$. 
%
\\

\cite[Proposition 4.3]{KauffmanNeumann} gives an answer to the question 
which we posed several paragraphs above.  
It shows that $L(f)\otimes L(g)$ is isotopic to $L(f + g)$.   
\\

Note that knot products $J\otimes K$ are defined 
even when neither $J$ nor $K$ is a link of singularity. 
This fact is significant because it utilizes structure that formerly was only available through algebraic varieties to knot theory proper.

\bigbreak
\section{Review of local moves on high dimensional knots}\label{move}

\begin{figure}
\unitlength 0.1in
\begin{picture}(36.39,15.20)(4.01,-22.91)
%
\special{pn 8}%
\special{ar 1138 1552 737 739  0.0000000 6.2831853}%
%
\special{pn 8}%
\special{pa 566 2017}%
\special{pa 1499 921}%
\special{fp}%
\special{sh 1}%
\special{pa 1499 921}%
\special{pa 1441 959}%
\special{pa 1464 962}%
\special{pa 1471 985}%
\special{pa 1499 921}%
\special{fp}%
%
\special{pn 8}%
\special{pa 1705 1085}%
\special{pa 771 2187}%
\special{fp}%
\special{sh 1}%
\special{pa 771 2187}%
\special{pa 829 2149}%
\special{pa 805 2146}%
\special{pa 799 2123}%
\special{pa 771 2187}%
\special{fp}%
%
\special{pn 8}%
\special{pa 1292 1510}%
\special{pa 1143 1401}%
\special{fp}%
%
\special{pn 8}%
\special{pa 1132 1707}%
\special{pa 988 1582}%
\special{fp}%
%
\special{pn 8}%
\special{pa 2375 1660}%
\special{pa 2091 1661}%
\special{fp}%
\special{sh 1}%
\special{pa 2091 1661}%
\special{pa 2158 1681}%
\special{pa 2144 1661}%
\special{pa 2158 1641}%
\special{pa 2091 1661}%
\special{fp}%
%
\special{pn 8}%
\special{pa 2102 1541}%
\special{pa 2401 1541}%
\special{fp}%
\special{sh 1}%
\special{pa 2401 1541}%
\special{pa 2334 1521}%
\special{pa 2348 1541}%
\special{pa 2334 1561}%
\special{pa 2401 1541}%
\special{fp}%
%
\put(18.7000,-21.7000){\makebox(0,0)[lb]{}}%
%
\special{pn 8}%
\special{ar 3303 1510 737 739  0.0000000 6.2831853}%
%
\special{pn 8}%
\special{pa 3148 1504}%
\special{pa 3272 1360}%
\special{fp}%
%
\special{pn 8}%
\special{pa 3298 1733}%
\special{pa 2968 2167}%
\special{fp}%
\special{sh 1}%
\special{pa 2968 2167}%
\special{pa 3024 2126}%
\special{pa 3000 2125}%
\special{pa 2992 2102}%
\special{pa 2968 2167}%
\special{fp}%
%
\special{pn 8}%
\special{pa 3102 1582}%
\special{pa 2751 1991}%
\special{fp}%
%
\special{pn 8}%
\special{pa 3339 1288}%
\special{pa 3664 874}%
\special{fp}%
\special{sh 1}%
\special{pa 3664 874}%
\special{pa 3607 914}%
\special{pa 3631 916}%
\special{pa 3639 939}%
\special{pa 3664 874}%
\special{fp}%
%
\special{pn 8}%
\special{pa 3344 1681}%
\special{pa 3489 1505}%
\special{fp}%
%
\special{pn 8}%
\special{pa 3535 1464}%
\special{pa 3865 1045}%
\special{fp}%
%
\special{pn 8}%
\special{pa 1365 1572}%
\special{pa 1787 1897}%
\special{fp}%
\special{sh 1}%
\special{pa 1787 1897}%
\special{pa 1746 1840}%
\special{pa 1745 1864}%
\special{pa 1722 1872}%
\special{pa 1787 1897}%
\special{fp}%
%
\special{pn 8}%
\special{pa 1087 1339}%
\special{pa 669 997}%
\special{fp}%
%
\special{pn 8}%
\special{pa 1632 2094}%
\special{pa 1184 1753}%
\special{fp}%
%
\special{pn 8}%
\special{pa 942 1536}%
\special{pa 503 1173}%
\special{fp}%
\special{sh 1}%
\special{pa 503 1173}%
\special{pa 542 1231}%
\special{pa 544 1207}%
\special{pa 567 1200}%
\special{pa 503 1173}%
\special{fp}%
%
\special{pn 8}%
\special{pa 2870 930}%
\special{pa 3958 1862}%
\special{fp}%
\special{sh 1}%
\special{pa 3958 1862}%
\special{pa 3920 1803}%
\special{pa 3917 1827}%
\special{pa 3894 1834}%
\special{pa 3958 1862}%
\special{fp}%
%
\special{pn 8}%
\special{pa 3772 2073}%
\special{pa 2664 1153}%
\special{fp}%
\special{sh 1}%
\special{pa 2664 1153}%
\special{pa 2703 1211}%
\special{pa 2705 1187}%
\special{pa 2728 1180}%
\special{pa 2664 1153}%
\special{fp}%
\end{picture}%

\caption{{\bf  The pass-move on 1-links}\label{Alabama}}
\end{figure}
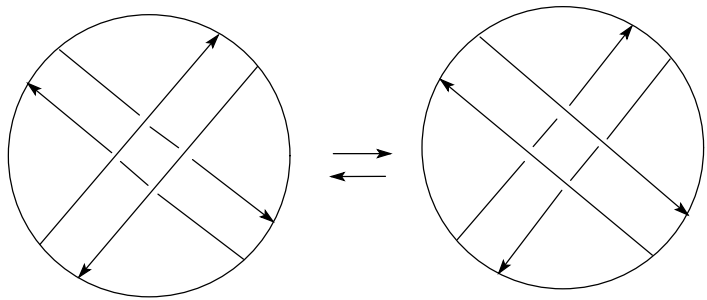

\noindent
We first review the definition of the pass-move on 1-links. 

\begin{defn}\label{Florida}   
({\bf\cite{Kauffmanon}.}) 
Let $K$ and $J$  be 1-links. 
If   $K$ (respectively,  $J$) is isotopic to a 1-link $K'$ (respectively,  $J'$) 
and 
if $K'$ and $J'$ differ only in a 3-ball  $B$ as shown in 
Figure \ref{Alabama},  
we say that $K$ (respectively,  $J$) is obtained from $J$ (respectively,  $K$) 
by one {\it pass-move}.  
If $K$ is obtained from a 1-link $P$ by a sequence of pass-moves, 
we say that $K$ is {\it pass-equivalent} to $P$. 
Each of four arcs in the 3-ball 
may belong to different components of the 1-link. 
\end{defn}

We use the terms `handle' and `surgeries' in this paper. 
See \cite{Browder, Kirby, Luck, Ranichi, Smale, Wall} 
for the definition of handles (respectively,  surgeries, the attaching parts of handles, 
the attached part, 
other related terms to handles).  
Note that 
an $a$-dimensional $q$-handle $h^q$ is diffeomorphic to $B^q\x B^{a-q}$ (respectively,  $B^a$), 
where $B^r$ denotes the $r$-ball, 
and that 
the attaching part of $h^q$ is diffeomorphic $S^{q-1}\x B^{a-q}$. 
Here, we review `a surgery by using an embedded handle'. 
In \S\ref{handleproof} we will review a few other facts on handle decompositions.

\begin{defn}\label{subsur}
Let $x$ and $m$ be positive integers 
and $x\leqq m$. 
Let $X$ be an $x$-dimensional submanifold of an $m$-dimensional manifold $M$.

Take $X\x [0,1]$ as an abstract manifold. 
Let $p$ be a nonnegative integer 
and $0\leqq p\leqq x$.
Attach an $(x+1)$-dimensional handle $h^p$ to $X\x [0,1]$. 
Suppose that the attaching part of $h^p$ is embedded only in $X\x\{1\}$. 
See Figure \ref{bat}.   
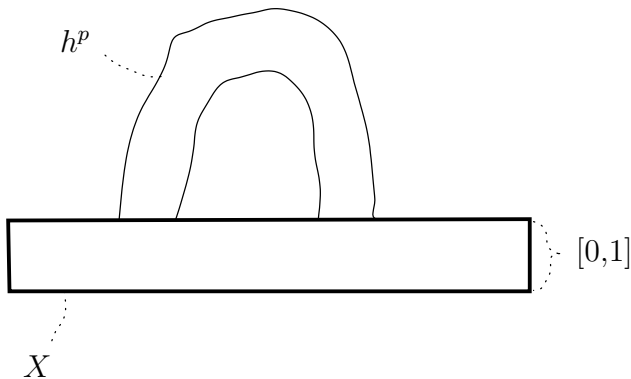
\begin{figure} 
\unitlength 0.1in
\begin{picture}(29.70,17.59)(16.60,-24.83)
%
\special{pn 20}%
\special{pa 1660 1834}%
\special{pa 4390 1827}%
\special{pa 4390 2205}%
\special{pa 1667 2205}%
\special{pa 1667 2205}%
\special{pa 1660 1834}%
\special{fp}%
%
\special{pn 8}%
\special{pa 2241 1827}%
\special{pa 2244 1795}%
\special{pa 2247 1763}%
\special{pa 2250 1731}%
\special{pa 2253 1700}%
\special{pa 2256 1668}%
\special{pa 2260 1636}%
\special{pa 2264 1604}%
\special{pa 2268 1572}%
\special{pa 2273 1541}%
\special{pa 2278 1509}%
\special{pa 2284 1477}%
\special{pa 2291 1446}%
\special{pa 2298 1415}%
\special{pa 2307 1384}%
\special{pa 2316 1353}%
\special{pa 2327 1323}%
\special{pa 2338 1293}%
\special{pa 2351 1264}%
\special{pa 2365 1235}%
\special{pa 2381 1208}%
\special{pa 2396 1181}%
\special{pa 2413 1154}%
\special{pa 2429 1127}%
\special{pa 2444 1100}%
\special{pa 2459 1071}%
\special{pa 2473 1042}%
\special{pa 2485 1012}%
\special{pa 2495 979}%
\special{pa 2503 945}%
\special{pa 2511 911}%
\special{pa 2523 882}%
\special{pa 2544 863}%
\special{pa 2574 853}%
\special{pa 2607 844}%
\special{pa 2637 831}%
\special{pa 2665 817}%
\special{pa 2693 802}%
\special{pa 2723 791}%
\special{pa 2753 783}%
\special{pa 2785 777}%
\special{pa 2817 771}%
\special{pa 2848 764}%
\special{pa 2879 757}%
\special{pa 2911 751}%
\special{pa 2942 746}%
\special{pa 2974 739}%
\special{pa 3005 731}%
\special{pa 3036 725}%
\special{pa 3068 724}%
\special{pa 3100 727}%
\special{pa 3132 733}%
\special{pa 3164 740}%
\special{pa 3195 748}%
\special{pa 3225 759}%
\special{pa 3255 771}%
\special{pa 3283 786}%
\special{pa 3310 804}%
\special{pa 3335 824}%
\special{pa 3359 847}%
\special{pa 3380 871}%
\special{pa 3399 897}%
\special{pa 3414 925}%
\special{pa 3427 954}%
\special{pa 3438 985}%
\special{pa 3447 1015}%
\special{pa 3455 1047}%
\special{pa 3463 1078}%
\special{pa 3470 1109}%
\special{pa 3479 1140}%
\special{pa 3488 1171}%
\special{pa 3497 1202}%
\special{pa 3505 1233}%
\special{pa 3514 1264}%
\special{pa 3522 1295}%
\special{pa 3530 1326}%
\special{pa 3537 1357}%
\special{pa 3543 1388}%
\special{pa 3548 1419}%
\special{pa 3553 1450}%
\special{pa 3556 1482}%
\special{pa 3560 1514}%
\special{pa 3563 1546}%
\special{pa 3566 1579}%
\special{pa 3569 1612}%
\special{pa 3572 1644}%
\special{pa 3575 1675}%
\special{pa 3578 1705}%
\special{pa 3580 1736}%
\special{pa 3572 1775}%
\special{pa 3567 1813}%
\special{pa 3585 1827}%
\special{pa 3585 1827}%
\special{sp}%
%
\special{pn 8}%
\special{pa 2535 1834}%
\special{pa 2545 1804}%
\special{pa 2555 1773}%
\special{pa 2564 1743}%
\special{pa 2573 1712}%
\special{pa 2582 1681}%
\special{pa 2591 1650}%
\special{pa 2599 1619}%
\special{pa 2607 1588}%
\special{pa 2614 1557}%
\special{pa 2620 1526}%
\special{pa 2625 1494}%
\special{pa 2630 1462}%
\special{pa 2634 1430}%
\special{pa 2636 1397}%
\special{pa 2640 1365}%
\special{pa 2646 1334}%
\special{pa 2656 1304}%
\special{pa 2670 1275}%
\special{pa 2686 1248}%
\special{pa 2703 1220}%
\special{pa 2719 1193}%
\special{pa 2736 1165}%
\special{pa 2755 1139}%
\special{pa 2778 1116}%
\special{pa 2805 1101}%
\special{pa 2836 1092}%
\special{pa 2868 1086}%
\special{pa 2899 1077}%
\special{pa 2930 1068}%
\special{pa 2960 1058}%
\special{pa 2992 1052}%
\special{pa 3024 1050}%
\special{pa 3056 1053}%
\special{pa 3086 1063}%
\special{pa 3114 1080}%
\special{pa 3138 1101}%
\special{pa 3161 1123}%
\special{pa 3183 1147}%
\special{pa 3202 1172}%
\special{pa 3219 1200}%
\special{pa 3232 1229}%
\special{pa 3243 1260}%
\special{pa 3251 1291}%
\special{pa 3256 1323}%
\special{pa 3259 1355}%
\special{pa 3260 1387}%
\special{pa 3260 1419}%
\special{pa 3261 1451}%
\special{pa 3264 1483}%
\special{pa 3268 1515}%
\special{pa 3274 1546}%
\special{pa 3281 1578}%
\special{pa 3287 1609}%
\special{pa 3291 1641}%
\special{pa 3294 1672}%
\special{pa 3294 1704}%
\special{pa 3292 1736}%
\special{pa 3290 1768}%
\special{pa 3286 1800}%
\special{pa 3284 1820}%
\special{sp}%
%
\special{pn 8}%
\special{pa 1954 2205}%
\special{pa 1958 2237}%
\special{pa 1961 2268}%
\special{pa 1961 2301}%
\special{pa 1956 2333}%
\special{pa 1945 2363}%
\special{pa 1926 2389}%
\special{pa 1908 2415}%
\special{pa 1895 2445}%
\special{pa 1885 2475}%
\special{pa 1884 2478}%
\special{sp -0.045}%
\put(17.3700,-26.5300){\makebox(0,0)[lb]{$X$}}%
%
\special{pn 8}%
\special{pa 2460 1080}%
\special{pa 2428 1082}%
\special{pa 2396 1081}%
\special{pa 2365 1073}%
\special{pa 2335 1061}%
\special{pa 2306 1048}%
\special{pa 2276 1036}%
\special{pa 2246 1024}%
\special{pa 2218 1009}%
\special{pa 2192 990}%
\special{pa 2168 969}%
\special{pa 2160 960}%
\special{sp -0.045}%
\put(19.3000,-9.5000){\makebox(0,0)[lb]{$h^p$}}%
%
\special{pn 8}%
\special{pa 4410 1840}%
\special{pa 4440 1850}%
\special{pa 4464 1871}%
\special{pa 4478 1900}%
\special{pa 4479 1932}%
\special{pa 4481 1964}%
\special{pa 4500 1990}%
\special{pa 4527 2005}%
\special{pa 4561 2010}%
\special{pa 4562 2007}%
\special{pa 4516 2008}%
\special{pa 4497 2032}%
\special{pa 4499 2065}%
\special{pa 4500 2097}%
\special{pa 4490 2125}%
\special{pa 4473 2154}%
\special{pa 4457 2186}%
\special{pa 4432 2200}%
\special{pa 4410 2200}%
\special{sp -0.045}%
\put(46.3000,-20.9000){\makebox(0,0)[lb]{[0,1]}}%
\end{picture}%
\bigbreak\bigbreak\quad\quad
\caption{   
{\bf A handle $h^p$ is attached to $X\x[0,1]$.}   
\label{bat}}
\bigbreak\end{figure}
\noindent 
Let 
$X'$\\
=$\overline{\partial(h^p\cup(X\x[0,1]))-(X\x\{0\})-((\partial X)\x[0,1])}$.
Note that there are two cases, $\partial X\\=\phi$ and $\partial X\neq\phi$. 

Suppose that there is a continuous injective map  $f:(X\x\{1\})\cup h^p\to M$ with the following properties: 
$f|_{X\x\{1\}}$, $f|_{h^p}$ and $f|_{X'}$ are smooth embedding maps. 
We have $f(X\x\{1\})=X$. Call the submanifold $f(X')$, $X'$ again. 

Then we say that 
the submanifold $X'$ is obtained from the submanifold $X$ 
by {\it the surgery by using the embedded handle $h^p$}. 
\end{defn}

\begin{figure} 
\includegraphics[width=10cm]{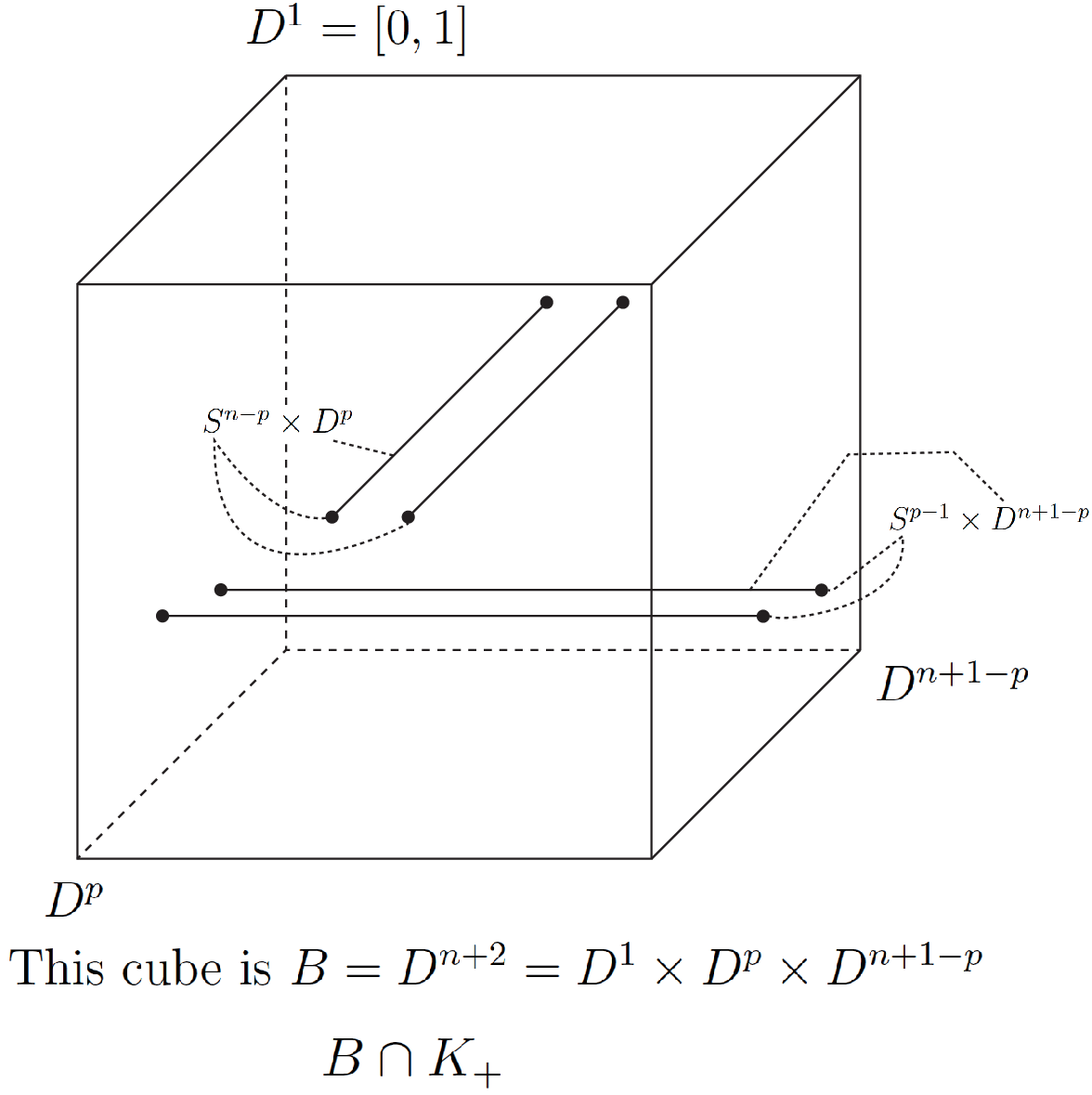}

\caption{{\bf  A $(p,n+1-p)$-pass-move triple}\label{Alaska}}
\end{figure}

\begin{figure} 
\includegraphics[width=10cm]{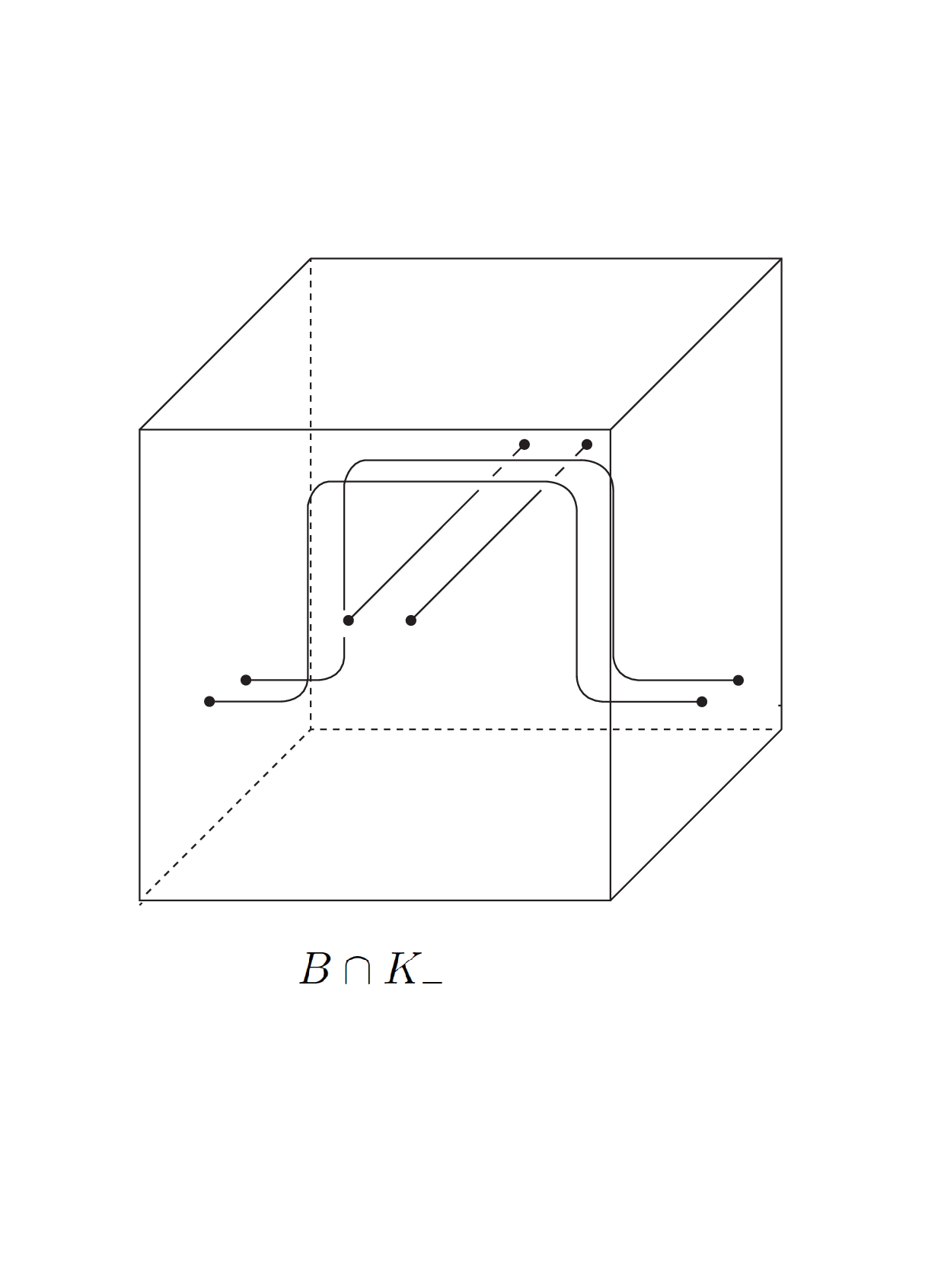}
\vskip-30mm
\caption{{\bf A $(p,n+1-p)$-pass-move triple}\label{Arizona}} 
\end{figure}

\begin{figure} 
\includegraphics[width=10cm]{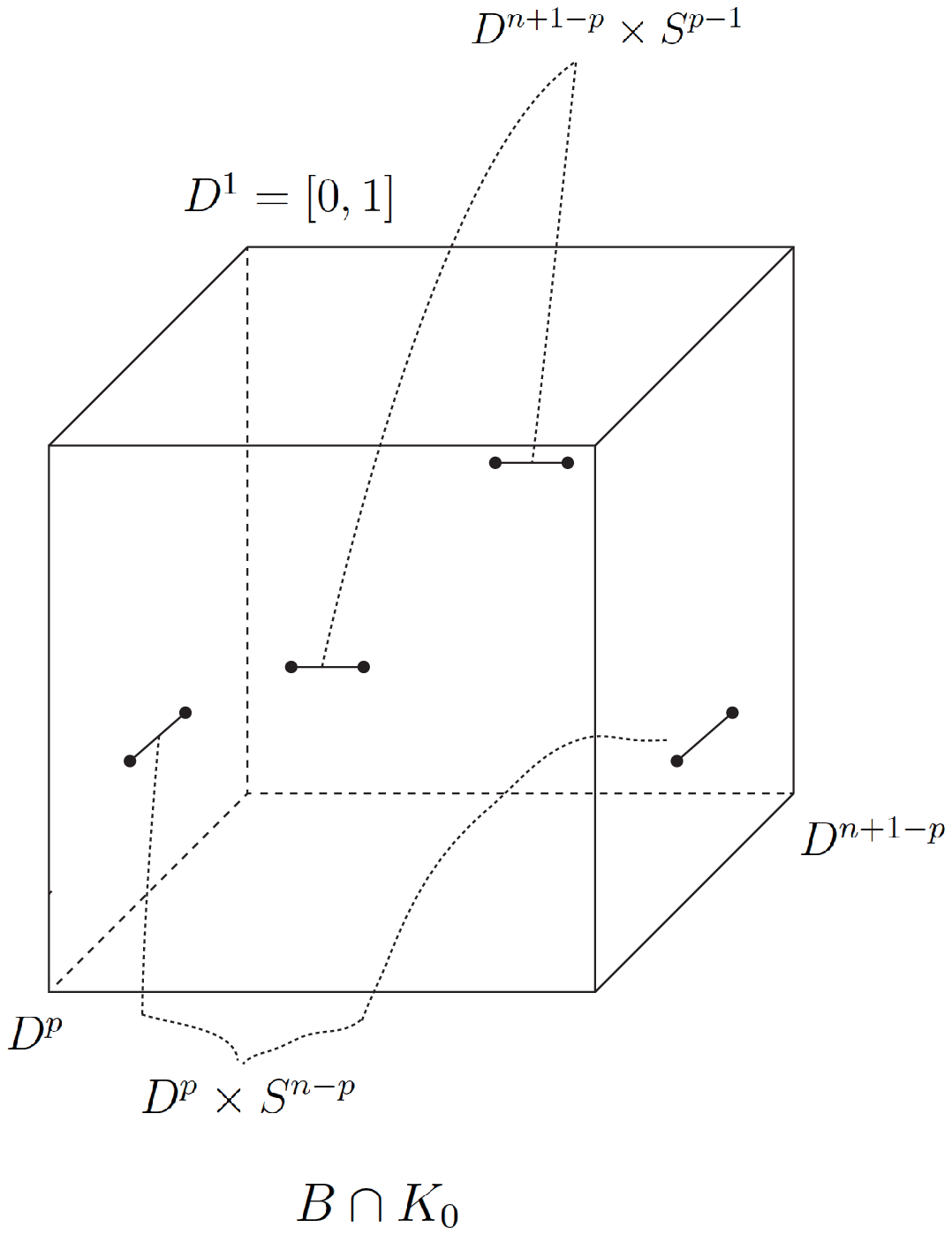}
\vskip-20mm
\caption{{\bf A $(p,n+1-p)$-pass-move triple}\label{Arkansas}}
\end{figure}

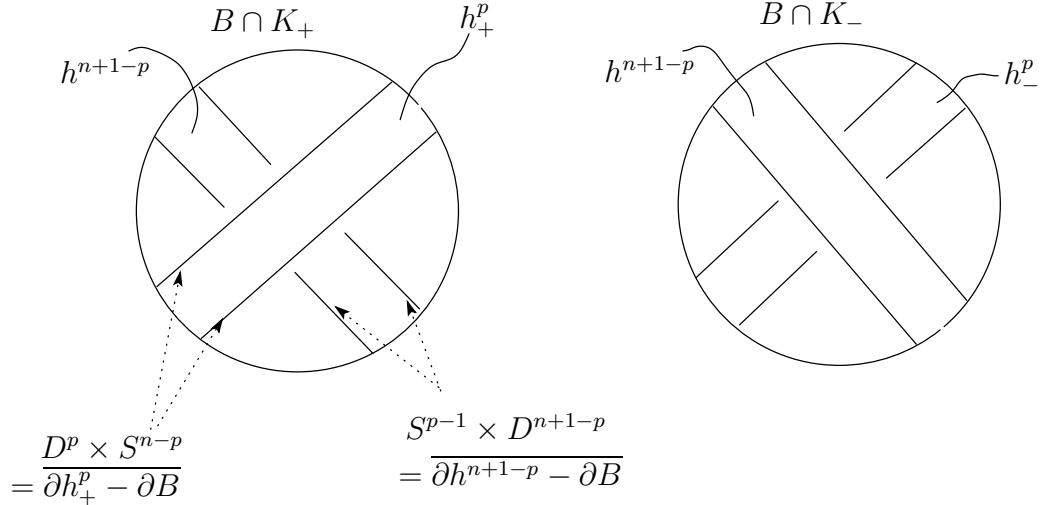
\begin{figure}
\unitlength 0.1in
\begin{picture}(52.00,24.84)(1.50,-32.30)
%
\special{pn 8}%
\special{ar 4488 1824 843 843  0.8934700 6.2831853}%
\special{ar 4488 1824 843 843  0.0000000 0.8617842}%
%
\special{pn 8}%
\special{pa 3829 1312}%
\special{pa 4891 2556}%
\special{fp}%
%
\special{pn 8}%
\special{pa 5155 2326}%
\special{pa 4104 1083}%
\special{fp}%
%
\special{pn 8}%
\special{pa 4891 1083}%
\special{pa 4515 1440}%
\special{fp}%
\special{pa 4168 1806}%
\special{pa 3738 2208}%
\special{fp}%
\special{pa 5147 1320}%
\special{pa 4735 1686}%
\special{fp}%
\special{pa 4369 2071}%
\special{pa 3966 2455}%
\special{fp}%
%
\special{pn 8}%
\special{ar 1652 1858 843 843  5.5899622 6.2831853}%
\special{ar 1652 1858 843 843  0.0000000 5.5598920}%
%
\special{pn 8}%
\special{pa 904 1468}%
\special{pa 1268 1838}%
\special{fp}%
\special{pa 1639 2180}%
\special{pa 2047 2604}%
\special{fp}%
\special{pa 1138 1209}%
\special{pa 1512 1614}%
\special{fp}%
\special{pa 1901 1974}%
\special{pa 2291 2371}%
\special{fp}%
%
\special{pn 8}%
\special{pa 2146 1183}%
\special{pa 917 2254}%
\special{fp}%
%
\special{pn 8}%
\special{pa 1149 2526}%
\special{pa 2377 1445}%
\special{fp}%
%
\special{pn 8}%
\special{pa 884 3002}%
\special{pa 1039 2179}%
\special{dt 0.045}%
\special{sh 1}%
\special{pa 1039 2179}%
\special{pa 1007 2241}%
\special{pa 1029 2231}%
\special{pa 1046 2248}%
\special{pa 1039 2179}%
\special{fp}%
\special{pa 921 2983}%
\special{pa 1258 2426}%
\special{dt 0.045}%
\special{sh 1}%
\special{pa 1258 2426}%
\special{pa 1206 2473}%
\special{pa 1230 2472}%
\special{pa 1241 2493}%
\special{pa 1258 2426}%
\special{fp}%
%
\special{pn 8}%
\special{pa 2402 2801}%
\special{pa 2228 2334}%
\special{dt 0.045}%
\special{sh 1}%
\special{pa 2228 2334}%
\special{pa 2233 2403}%
\special{pa 2247 2384}%
\special{pa 2270 2389}%
\special{pa 2228 2334}%
\special{fp}%
\special{pa 2393 2783}%
\special{pa 1853 2361}%
\special{dt 0.045}%
\special{sh 1}%
\special{pa 1853 2361}%
\special{pa 1893 2418}%
\special{pa 1895 2394}%
\special{pa 1918 2386}%
\special{pa 1853 2361}%
\special{fp}%
\put(22.1900,-30.6600){\makebox(0,0)[lb]{$S^{p-1}\x D^{n+1-p}$}}%
\put(3.0700,-31.6600){\makebox(0,0)[lb]{$D^p\x S^{n-p}$}}%
\put(21.8200,-33.0500){\makebox(0,0)[lb]{$=\overline{\partial h^{n+1-p}-\partial B}$}}%
\put(11.9400,-9.5300){\makebox(0,0)[lb]{$B\cap K_+$}}%
\put(40.6700,-9.1600){\makebox(0,0)[lb]{$B\cap K_-$}}%
\put(1.5000,-34.0000){\makebox(0,0)[lb]{$=\overline{\partial h^p_+-\partial B}$}}%
%
\special{pn 8}%
\special{pa 1100 1450}%
\special{pa 1111 1419}%
\special{pa 1121 1389}%
\special{pa 1128 1358}%
\special{pa 1130 1326}%
\special{pa 1128 1295}%
\special{pa 1122 1263}%
\special{pa 1112 1232}%
\special{pa 1100 1201}%
\special{pa 1087 1171}%
\special{pa 1071 1143}%
\special{pa 1053 1116}%
\special{pa 1033 1092}%
\special{pa 1009 1072}%
\special{pa 982 1055}%
\special{pa 954 1040}%
\special{pa 924 1026}%
\special{pa 893 1015}%
\special{pa 862 1004}%
\special{pa 831 1000}%
\special{pa 800 1007}%
\special{pa 770 1020}%
\special{sp}%
\put(4.1000,-11.9000){\makebox(0,0)[lb]{$h^{n+1-p}$}}%
\put(32.6000,-11.8000){\makebox(0,0)[lb]{$h^{n+1-p}$}}%
%
\special{pn 8}%
\special{pa 4040 1360}%
\special{pa 4035 1327}%
\special{pa 4028 1296}%
\special{pa 4018 1266}%
\special{pa 4002 1238}%
\special{pa 3982 1214}%
\special{pa 3961 1189}%
\special{pa 3941 1165}%
\special{pa 3921 1140}%
\special{pa 3902 1115}%
\special{pa 3883 1089}%
\special{pa 3864 1062}%
\special{pa 3844 1035}%
\special{pa 3821 1012}%
\special{pa 3794 998}%
\special{pa 3761 994}%
\special{pa 3729 1000}%
\special{pa 3700 1015}%
\special{pa 3670 1020}%
\special{sp}%
%
\special{pn 8}%
\special{pa 2190 1380}%
\special{pa 2199 1349}%
\special{pa 2209 1318}%
\special{pa 2219 1288}%
\special{pa 2231 1259}%
\special{pa 2245 1230}%
\special{pa 2261 1202}%
\special{pa 2280 1176}%
\special{pa 2300 1151}%
\special{pa 2322 1127}%
\special{pa 2346 1106}%
\special{pa 2373 1089}%
\special{pa 2404 1078}%
\special{pa 2435 1069}%
\special{pa 2464 1058}%
\special{pa 2488 1040}%
\special{pa 2508 1016}%
\special{pa 2524 987}%
\special{pa 2539 956}%
\special{pa 2550 930}%
\special{sp}%
\put(25.0000,-9.5000){\makebox(0,0)[lb]{$h^p_+$}}%
\put(53.5000,-12.4000){\makebox(0,0)[lb]{$h^p_-$}}%
%
\special{pn 8}%
\special{pa 5000 1290}%
\special{pa 5028 1274}%
\special{pa 5054 1255}%
\special{pa 5078 1234}%
\special{pa 5098 1208}%
\special{pa 5117 1180}%
\special{pa 5140 1160}%
\special{pa 5170 1152}%
\special{pa 5203 1145}%
\special{pa 5234 1140}%
\special{pa 5265 1145}%
\special{pa 5297 1150}%
\special{pa 5320 1150}%
\special{sp}%
\end{picture}%
\caption{{\bf  The $(p, n+1-p)$-move on an $n$-dimensional closed submanifold of $S^{n+2}$. 
Note $B=B^{n+2}=D^{n+2}\subset S^{n+2}$.}\label{Georgia}}
\bigbreak
\end{figure}

Let $p, q\in\N$ and $p+q=n+1$. 
We review the definition of the $(p, q)$-pass-move on $n$-submanifolds of $S^{n+2}$, 
which was defined in 
\cite{Ogasa98n} 
and 
which has been studied in 
\cite{KauffmanOgasa, Ogasa98n,  Ogasa02,  Ogasa04, Ogasa07,   Ogasa09, OgasaT3, OgasaIH}. 
The triple of Figures 
\ref{Alaska}-\ref{Arkansas}   
is a diagram 
of a $(p, q)$-pass-move triple. 
Figure \ref{Georgia} also 
represents a diagram of the $(p, q)$-pass-move on $n$-submanifolds of $S^{n+2}$,  where $q=n+1-p$.  
Note 
that, if $(p, q)=(1, 1)$, 
the $(p, q)$-pass-move is the pass-move on 1-links in Figure \ref{Alabama}. 
If $p=1$ and $q=2$, any $(p,q)$-pass-move triple of $(p+q-1)$-submanifolds of $S^{p+q+1}$ is  
a $(1,2)$-pass-move triple of 2-submanifolds of $S^4$, which is defined in \cite{Ogasa09}. 
%
%
%
%
In Definition \ref{Minesota}  we explain the definition of 
a $(p,q)$-pass-move triple in more detail.

\begin{defn}\label{Minesota}{\bf(\cite{Ogasa98n}.)}
Let $n,p\in\N$. Let $n+1-p>0.$     
We now define the $(p, n+1-p)$-pass-move in a $(n+2)$-ball. We first explain 
Figures \ref{Alaska}-\ref{Arkansas}.   \\

Regard an $(n+2)$-ball $B=D^{n+2}$ as $[-1,1]\x D^p\x D^{n+1-p}$ 
as drawn in 
Figures \ref{Alaska}. \\

Attach an embedded $(n+1)$-dimensional $(n+1-p)$-handle 
$h^{n+1-p} \subset\{0\}\x D^p\x D^{n+1-p}$ 
to 
$\{0\}\x\partial(D^p\x D^{n+1-p})$
along 
$\{0\}\x\{*\}\x \partial D^{n+1-p}$
so that 
the core of $h^{n+1-p}$ coincides with 
$\{0\}\x\{*\}\x D^{n+1-p}$.  
See Figures \ref{Alaska} and \ref{Arizona}. \\ 

Attach an embedded $(n+1)$-dimensional $p$-handle 
$h^p \subset\{0\}\x D^p\x D^{n+1-p}$ 
to \newline
$\{0\}\x\partial(D^p\x D^{n+1-p})$
along 
$\{0\}\x \partial D^p\x\{*\}$
so that 
the core of $h^p$ coincides with 
$\{0\}\x D^p\x\{*\}$.  
Note $h^{n+1-p}\cap h^p\neq\phi.$\\

Move $h^p$ in $B$ by using an isotopy with keeping $h^p\cap \partial B$, 
 let 
the resultant submanifold of $\{t\geqq0\}\x D^p\x D^{n+1-p}$ 
(respectively,  $\{t\leqq0\}\x D^p\x D^{n+1-p}$), 
and call the submanifold, $h_+^p$ (respectively,  $h_-^p$).  \\

Suppose that $h_+^p\cap h^{n+1-p}=\phi$ and 
that $h_-^p\cap h^{n+1-p}=\phi$.   
See Figures \ref{Alaska} and \ref{Arizona}. \\
Call $h^p\cap \partial B=h_+^p\cap \partial B=h_-^p\cap \partial B$, $P$. 
Call $h^{n+1-p}\cap \partial B$, $Q$. \\

Let $K_+$, $K_-$, and $K_0$ be $n$-dimensional closed oriented submanifolds of $S^{n+2}$. 
Embed the $(n+2)$-ball $B$ in $S^{n+2}$. 
Let $K_+$, $K_-$ and $K_0$ differ only in $B$. 
Let $\amalg$ denote the disjoint union. 
Let $K_+$ (respectively,  $K_-$, $K_0$) satisfy the condition 
$$K_+\cap \mathrm{Int}B= 
(\partial h^{p}_+-\partial B)\amalg(\partial h^{n+1-p}-\partial B)$$
$${\rm(respectively, }\hskip1mm K_-\cap \mathrm{Int}B=
(\partial h^{p}_--\partial B)\amalg(\partial h^{n+1-p}-\partial B), 
K_0\cap B=P\amalg Q)$$

\noindent
where we suppose that there is not  
$h^{p}_-$ or $h^{p}$  
(respectively,  `$h^{p}_+$ or $h^{p}$',  `$h^{n+1-p}$, $h^{p}_+$, or $h^{p}$')   
in $B$.  
Then we say that 
$K_+$ (respectively,  $K_-$)  is obtained from $K_-$ (respectively,  $K_+$) 
by one {\it $(p,n+1-p)$-pass-move} in $B$. 
We say that  $(K_+$, $K_-$, $K_0)$ is related by a single {\it $(p,n+1-p)$-pass-move} in $B$.   
We also say that $(K_+$, $K_-$, $K_0)$ is a {\it $(p,n+1-p)$-pass-move triple}.   
If \newline $(K_+$, $K_-$, $K_0)$ is a $(p,n+1-p)$-pass-move triple, 
then we also say that $(K_-$, $K_+$, $K_0)$ is a {\it $(p,n+1-p)$-pass-move triple}.   
\end{defn}

In Definition \ref{Minesota}, we have the following: 
Let $\sharp\in\{+,-\}.$   
There is a Seifert hypersurface $V_\sharp\subset S^{n+2}$ for $K_\sharp$ such that 
$V_\sharp\cap B=h_\sharp^p\cup h^{n+1-p}.$ 
\noindent 
(The idea of the proof is Thom-Pontrjagin construction.)
We say that  
$V_-$ (respectively,  $V_+$) is obtained from $V_+$ (respectively,  $V_-$) 
by a {\it $(p,n+1-p)$-pass-move} in $B$.  

Let 
$V_0=V_\sharp-\text{Int}B$
$=\text{the closure of }`V_\sharp- (h_\sharp^p\cup{h^{n+1-p}})'\text{ in }{S^{n+2}}.$  
We can say that 
we attach 
an embedded $(n+1)$-dimensional $p$-handle $h_\#^p\subset S^{n+2}$ 
and 
an embedded $(n+1)$-dimensional $(n+1-p)$-handle $h^{n+1-p}\subset S^{n+2}$ 
to the $(n+1)$-submanifold $V_0\subset S^{n+2}$, 
and 
obtain the $(n+1)$-submanifold $V_\#\subset S^{n+2}$.





The ordered set ($V_+, V_-, V_0$) is called a {\it $(p,n+1-p)$-pass-move triple of 
Seifert hypersurfaces} for  $(K_+, K_-, K_0)$. 
We also say that 
an ordered set ($V_+, V_-, V_0$) is related by a single {\it $(p,n+1-p)$-pass-move} in $B$. 


\begin{figure}
\includegraphics[width=9cm]{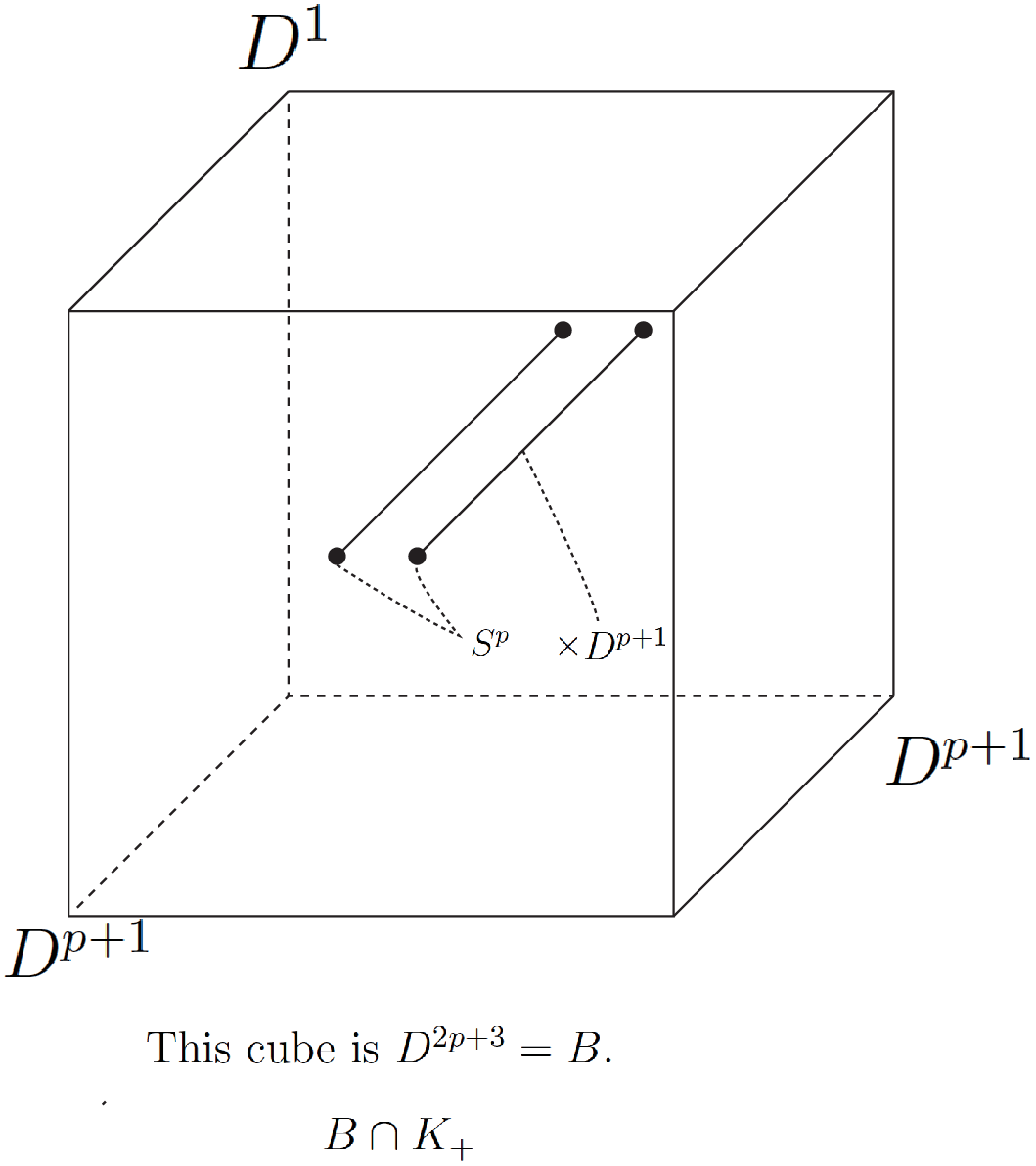}
\caption{{\bf  A twist-move triple}\label{California}} 
\end{figure}

\begin{figure}
\includegraphics[width=9cm]{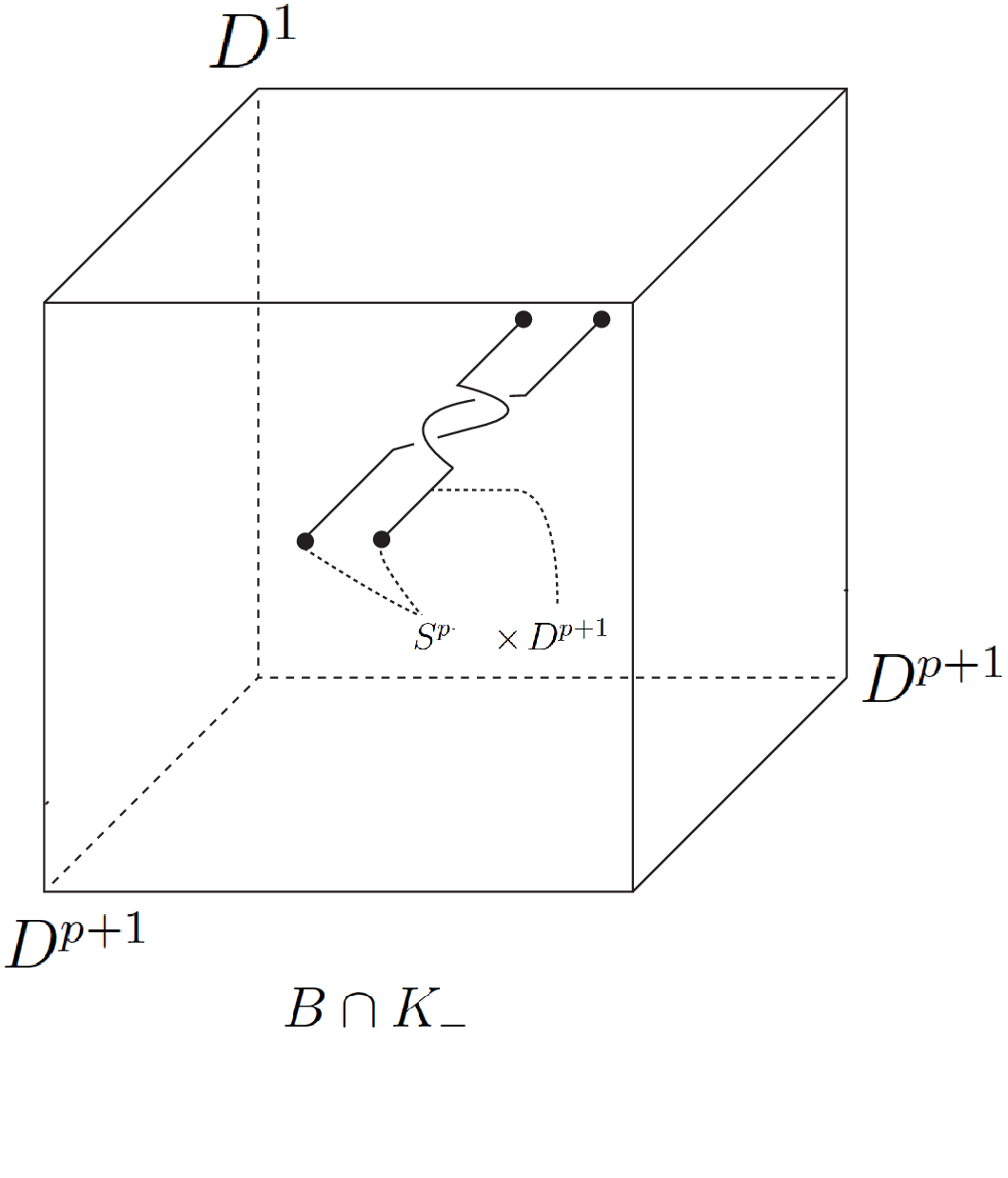}
\vskip-10mm
\caption{{\bf A twist-move triple}\label{Colorado}} 
\end{figure}

\begin{figure}
\includegraphics[width=10cm]{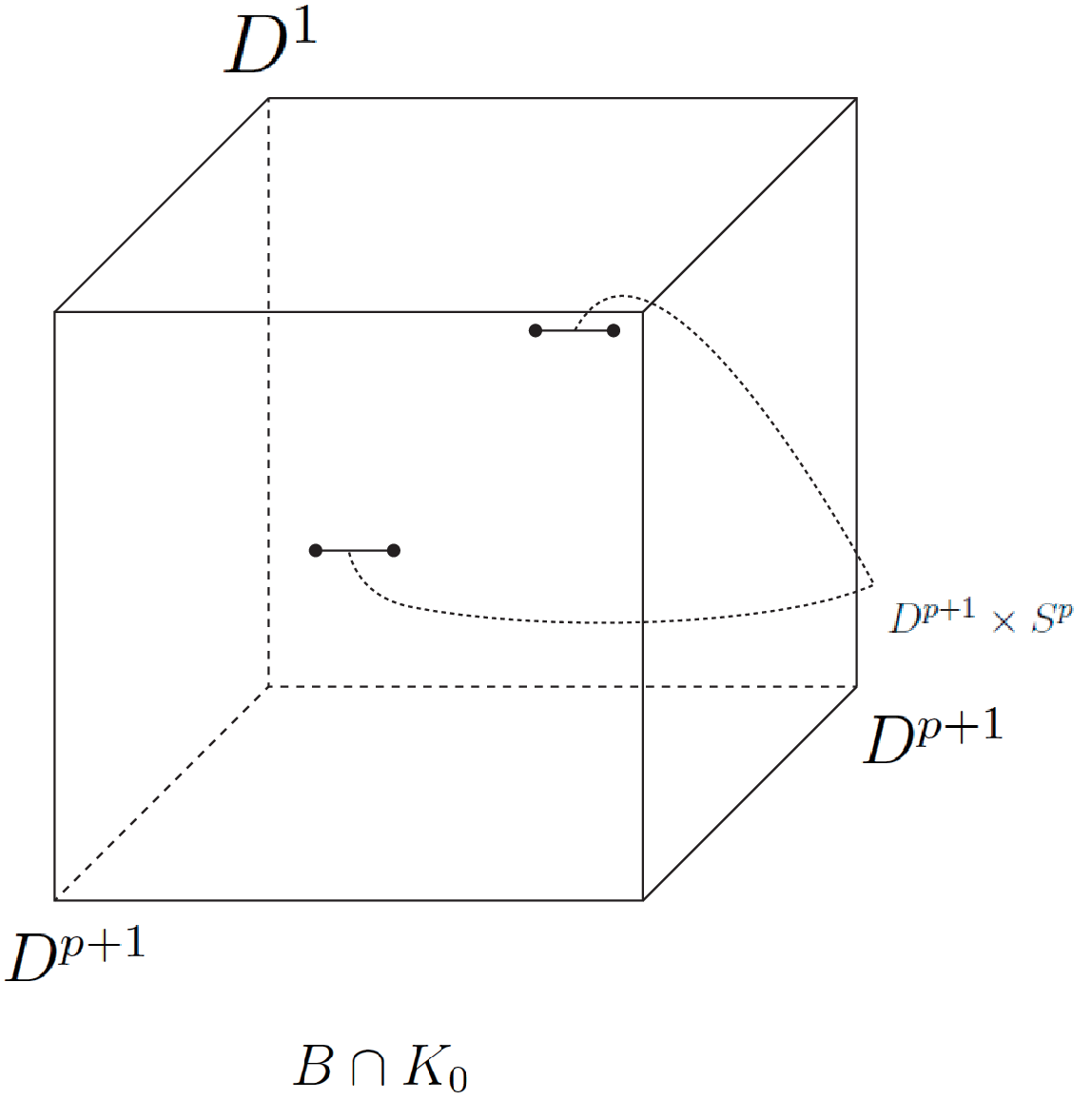}
\caption{{\bf A twist-move triple}\label{Connecticut}} 
\end{figure}

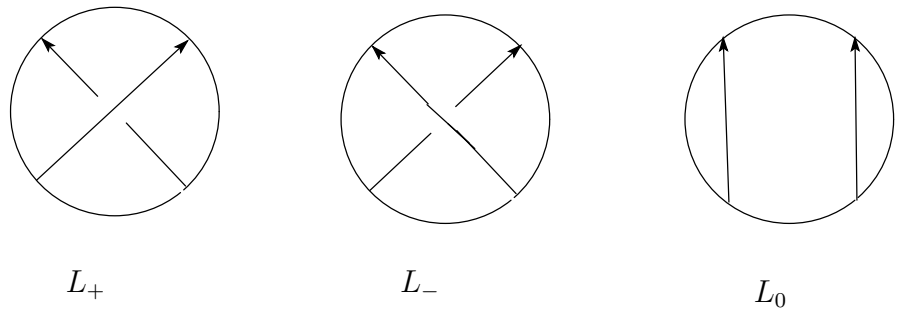
\begin{figure}
\unitlength 0.1in
\begin{picture}(46.22,14.00)(8.10,-22.10)
%
\special{pn 8}%
\special{ar 1356 1356 546 546  0.8902751 6.2831853}%
\special{ar 1356 1356 546 546  0.0000000 0.8502422}%
%
\special{pn 8}%
\special{pa 946 1716}%
\special{pa 1736 986}%
\special{fp}%
\special{sh 1}%
\special{pa 1736 986}%
\special{pa 1673 1017}%
\special{pa 1697 1022}%
\special{pa 1701 1046}%
\special{pa 1736 986}%
\special{fp}%
%
\special{pn 8}%
\special{pa 1266 1276}%
\special{pa 976 976}%
\special{fp}%
\special{sh 1}%
\special{pa 976 976}%
\special{pa 1008 1038}%
\special{pa 1013 1014}%
\special{pa 1037 1010}%
\special{pa 976 976}%
\special{fp}%
%
\special{pn 8}%
\special{pa 1416 1416}%
\special{pa 1726 1746}%
\special{fp}%
%
\special{pn 8}%
\special{ar 3086 1396 546 546  0.8902751 6.2831853}%
\special{ar 3086 1396 546 546  0.0000000 0.8502422}%
%
\special{pn 8}%
\special{pa 2996 1316}%
\special{pa 2706 1016}%
\special{fp}%
\special{sh 1}%
\special{pa 2706 1016}%
\special{pa 2738 1078}%
\special{pa 2743 1054}%
\special{pa 2767 1050}%
\special{pa 2706 1016}%
\special{fp}%
%
\special{pn 8}%
\special{pa 3146 1456}%
\special{pa 3456 1786}%
\special{fp}%
%
\special{pn 8}%
\special{pa 3210 1520}%
\special{pa 2990 1320}%
\special{fp}%
\special{pa 3240 1550}%
\special{pa 3100 1420}%
\special{fp}%
%
\special{pn 8}%
\special{pa 3140 1330}%
\special{pa 3480 1010}%
\special{fp}%
\special{sh 1}%
\special{pa 3480 1010}%
\special{pa 3418 1041}%
\special{pa 3441 1047}%
\special{pa 3445 1070}%
\special{pa 3480 1010}%
\special{fp}%
%
\special{pn 8}%
\special{pa 3010 1470}%
\special{pa 2690 1770}%
\special{fp}%
%
\special{pn 8}%
\special{ar 4886 1396 546 546  0.8902751 6.2831853}%
\special{ar 4886 1396 546 546  0.0000000 0.8502422}%
%
\special{pn 8}%
\special{pa 4570 1830}%
\special{pa 4540 970}%
\special{fp}%
\special{sh 1}%
\special{pa 4540 970}%
\special{pa 4522 1037}%
\special{pa 4542 1023}%
\special{pa 4562 1036}%
\special{pa 4540 970}%
\special{fp}%
%
\special{pn 8}%
\special{pa 5240 1820}%
\special{pa 5230 970}%
\special{fp}%
\special{sh 1}%
\special{pa 5230 970}%
\special{pa 5211 1037}%
\special{pa 5231 1023}%
\special{pa 5251 1036}%
\special{pa 5230 970}%
\special{fp}%
\put(11.0000,-23.4000){\makebox(0,0)[lb]{$L_+$}}%
\put(28.5000,-23.4000){\makebox(0,0)[lb]{$L_{-}$}}%
\put(47.0000,-23.8000){\makebox(0,0)[lb]{$L_0$}}%
\end{picture}%
\caption{{\bf A crossing change triple of 1-links}\label{Hawaii}}   
\smallbreak  
\end{figure}

\begin{figure}
\includegraphics[width=14cm]{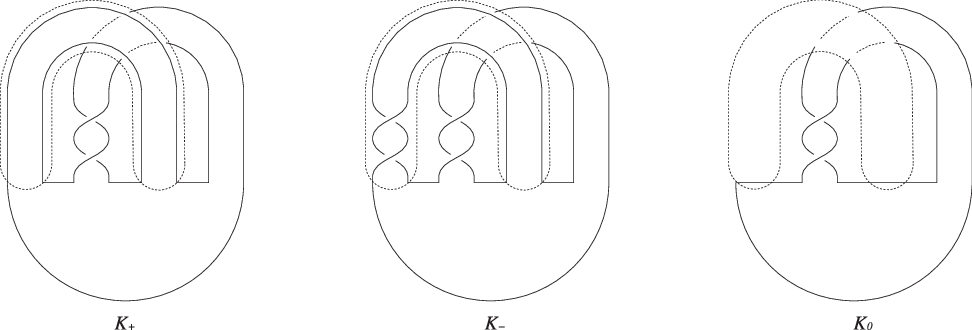}
\caption{{\bf A twist-move triple of $(4k+1)$-knots}\label{Oregon}}   
\end{figure}

\bigbreak
We review the definition of 
the twist-move on $(2p+1)$-submanifolds of $S^{2p+3}$,  
which is defined in \cite{Ogasa09}. 
Note that, there, the twist-move is called the $XXII$-move. 
The triple of Figures 
\ref{California}-\ref{Connecticut}   
is a diagram of a twist-move triple in a $(2p+3)$-ball in the standard $(2p+3)$-sphere. 
Figure \ref{Oregon} is an example of 
a twist-move triple of $(2p+1)$-submanifolds of $S^{2p+3}$. 
Note 
that if $p=0$, 
any twist-move triple of $(2p+1)$-submanifolds of $S^{2p+3}$  
is a crossing change triple of 1-links 
in Figure \ref{Hawaii}.     

\begin{defn}\label{0spoon}  
Let $p\in\N\cup\{0\}$.  
We now define the twist-move in a $(2p+3)$-ball. 
We first explain Figures \ref{California}-\ref{Connecticut}.    \\

Regard a $(2p+3)$-ball $B=D^{2p+3}$ as $[-1,1]\x D^{p+1}\x D^{p+1}$
as drawn in Figures \ref{California}-\ref{Connecticut}.   \\     

Attach an embedded $(2p+2)$-dimensional $(p+1)$-handle 
$h_+\subset\{0\}\x D^{p+1}\x D^{p+1}$ 
(respectively,  $h_-\subset[-1,1]\x D^{p+1}\x D^{p+1}$) 
to 
$\{0\}\x\partial(D^{p+1}\x D^{p+1})$
along 
$\{0\}\x\{*\}\x \partial D^{p+1}$
so that 
the core of $h_+$ (respectively,  $h_-$) coincides with 
$\{0\}\x\{*\}\x D^{p+1}$.  
See Figures \ref{California} and \ref{Colorado}.     
We can suppose that $h_+\cap h_-$ 
is the attaching part of $h_+$ (respectively,  $h_-$), and  call $h_+\cap h_-$, $Q$.   \\

We give an orientation to $h_+$.  
We give an orientation to $h_-$ so that 
 $h_+\cup h_-$ is an oriented submanifold of $B$  
if we give the opposite orientation to $h_-$.
We can regard $h_+\cup h_-$ as a Seifert hypersurface for 
a closed oriented $(2p+1)$-submanifold $\partial(h_+\cup h_-)\subset B=D^{2p+3}$. 
We can suppose that 
a $(p+1)$-Seifert matrix 
for a $(2p+1)$-dimensional closed oriented submanifold $\partial(h_+\cup h_-)\subset B$ 
associated with a Seifert hypersurface $h_+\cup h_-$ is $(1)$. 
(We can define Seifert hypersurfaces in $B$ and their Seifert matrices  
in the same fashion as ones in the $S^n$ case.)  \\

Let  $p\in\N\cup\{0\}$. 
Let $K_+$, $K_-$, and $K_0$ be 
$(2p+1)$-dimensional closed oriented submanifold of $S^{2p+3}$.
Take $B$ in $S^{2p+3}$. 
Let $K_+$, $K_-$ and $K_0$ differ only in $B$.
Let $K_+$ (respectively,  $K_-$, $K_0$) satisfy the condition  
$$K_+\cap \mathrm{Int}B=\partial h_+-\partial B$$
$$({\rm respectively,  } \hskip1mm K_-\cap \mathrm{Int}B=\partial h_--\partial B,  K_0\cap B=Q)$$

\noindent
where we suppose that there is not 
$h_-$   
(respectively,  $h_+$, `$h_+$ or $h_-$') in $B$.  
Then we say that 
$K_+$ (respectively,  $K_-$)  is obtained from $K_-$ (respectively,  $K_+$) 
by one 
{\it $($positive$)$ twist-move} 
(respectively,  {\it $($negative$)$ twist-move})  
in $B$.  
We say that  an ordered set $(K_+$, $K_-$, $K_0)$ is related by a single {\it twist-move} 
and that $(K_+$, $K_-$, $K_0)$ is a {\it twist-move triple}.
\end{defn}

Figure \ref{Oregon} is an example of a twist-move triple of $(4k+1)$-submanifolds of $S^{4k+3}$.\\

Let $(K_+, K_-, K_0)$ be related by a single twist-move in $B$.   
Then there is a Seifert hypersurface $V_*$ for $K_*$ ($*=+, -, 0$)  with the following properties. 

\smallbreak\noindent 
(1) $V_\sharp=V_0\cup h_\sharp$  ($\sharp=+,-$).
 $V_\sharp\cap B=h_\sharp$.  

\smallbreak
 \noindent 
(2) 
$V_0\cap$ Int $B=\phi$. 
%
$V_0\cap\partial B=Q$ 

\noindent 
(The idea of the proof is the Thom-Pontrjagin construction.)

\bigbreak
 The ordered set ($V_+, V_-, V_0$) is called a 
{\it twist-move triple of Seifert hypersurfaces} for  $(K_+, K_-, K_0)$. 
We say that 
$V_-$ (respectively,  $V_+$) is obtained from $V_+$ (respectively,  $V_-$) 
by a single 
 {\it negative-twist move} 
(respectively,  {\it positive-twist move})  in $B$. \\

\begin{figure}
\includegraphics[width=14cm]{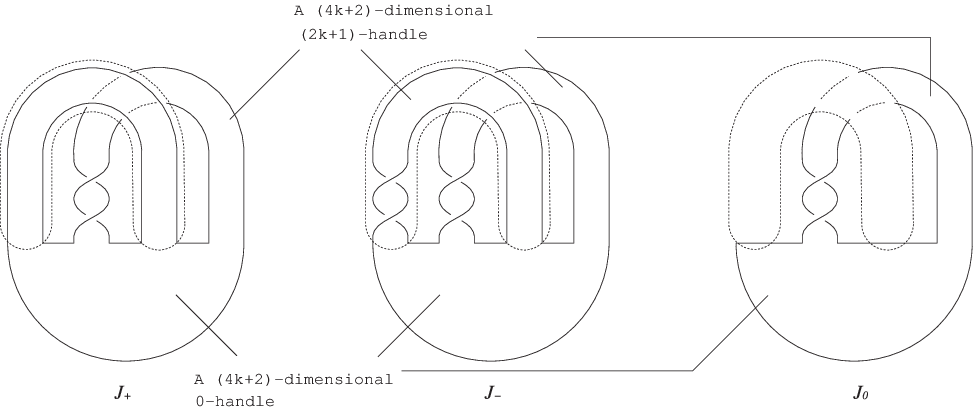}
\caption{{\bf A twist-move triple of Seifert hypersurfaces for \newline
$(4k+1)$-knots}\label{Oklahoma}} 
\end{figure}

Figure \ref{Oklahoma} is an example of a twist-move triple of Seifert hypersurfaces for $(4k+1)$-knots. \\

\smallbreak
\noindent{\bf{Note.}} 
Suppose that $p$ is an odd natural number, put $p=2k+1$. 
The twist-move for $(4k+3)$-submanifolds of $S^{4k+5}$ 
($k\in\N\cup\{0\}$)   
has the following property: 
Suppose that $K_+$ is made into $K_-$  by one twist-move. 
Then $K_-$ is a nonspherical knot in general even if $K_+$ is a spherical knot. 
Furthermore the   $H_*(K_-;\Z)$ is not congruent to $H_*(K_+;\Z)$ in general.   
Example: A Seifert hypersurface $V_*$ for  a 3-knot $K_*$ ($*=+,-$); 
Framed link representation of $V_+$  is the Hopf link 
such that the framing of one component is zero 
and that that of the other is two. 
Framed link representation of $V_-$  is the Hopf link 
such that the framing of each component is two.

\bigbreak
\section{Main results on local moves on knots and products of knots, 
and review of the authors' old results 
}\label{LP}  

We state our main results 
(Theorems \ref{Bunkyo} and \ref{aoiro}) on local moves on knots and products of knots. 
We review the authors' old results. We generalize the old results and obtain the main results. 
\\

\noindent
{\bf Knot products.}   
In \cite{Kauffman,  KauffmanNeumann} 
the following is proved: 
For any spherical $n$-knot $K$ in $S^{n+2}$ ($n\in\N$),    
take $K\otimes$Hopf, where we abbreviate $A\otimes$(the Hopf link) to $A\otimes$Hopf. 
Then we obtain a homomorphism 
$C_n\to C_{n+4}$, where $C_n$ is the knot cobordism group in \cite{Levinecob}.
Furthermore it is an isomorphism if $n\in\N-\{1,3\}$. 
We also obtained results in the case $n=1,3.$
\\

In \cite{KauffmanOgasa, KauffmanOgasaII} 
 we proved the following: 
For any simple $n$-knot $K$ in $S^{n+2}$ ($n\in\N$), take $K\otimes$Hopf. 
Then we obtain a map 
$\mathcal S_n\to \mathcal S_{n+4}$, where $\mathcal S_n$ is 
a set of simple $n$-knots.  
(See \S\ref{Introduction} and \ref{RkB}  for the definition of simple knots.)
Furthermore it is the bijective map 
if $n=5$ or $n\geqq7$.   
We also obtained results in the case $n=1,2,3,4,6.$
\\

In \cite{Kauffman,  KauffmanNeumann}  the empty knot $[a]$ of degree $a$ ($a\in\N$) was introduced. 
It is proved that 
any Brieskorn submanifold is isotopic to a $[a_1]\otimes...\otimes[a_p]$ and 
that any $[a_1]\otimes...\otimes[a_p]$ is isotopic to a Brieskorn submanifold.

\bigbreak  
\noindent
{\bf Local moves on high dimensional knots.}  
In \cite{Ogasa98n}  the following is proved: 
The $(p+1,p+1)$-pass-move on spherical $(2p+1)$-knots 
preserves the Arf invariant (respectively,  the signature) if $p$ is even (respectively,  odd). 
Furthermore the following is proved: Let $p$ be even (respectively,  odd). 
Simple $(2p+1)$-knots, $K$ and $J$, are $(p+1,p+1)$-pass-equivalent 
if and only if 
their Arf invariant (respectively,  their signature) are the same.

\cite[A result between Definitions 5.1 and 5.2]{Ogasaribbontwo} proved the following: 
Let $p$ be any even (respectively,  odd) positive integer. 
Let $K$ be a $(2p+1)$ knot.  
$K$ is $(p+1,p+1)$-pass-move equivalent to the unknot 
if and only if 
$K$ is a simple knot and 
the Arf invariant (respectively,  the signature) is zero.  
\\

In \cite{Ogasa04, Ogasa07,   Ogasa09} 
the following are proved: 
The (1,2)-pass move (respectively,  the ribbon-move) on spherical 2-knots preserves 
the $\mu$-invariant of 2-knots, 
the $\Q/\Z$-valued $\widetilde\eta$-invariants of 2-knots, 
the Farber-Levine pairing of 2-knots,   
and 
partial information of the cup product of three elements in 
$H^1_{cpt}($the complement of each 2-knot$)$.   
Note that the ribbon-move is a kind of local moves on 
2-submanifolds of $S^4$. 
\\

In \cite{OgasaT3} 
the following is proved:  
For the Alexander polynomial $A(\hskip2mm)$ of high dimensional knots 
we have an identity 
$$A(K_+)-A(K_-)=(t-1)A(K_0)$$ 
associated with the twist move on $(4k+1)$-dimensional knots, where $k\in\N\cup\{0\}$  
(respectively,  the $(p,q)$-pass-move, where $p\neq q$, on high dimensional knots).   
\\

In \cite{KauffmanOgasa} the following is proved: 
For the Alexander polynomial $A(\hskip2mm)$ of 
$(4k+3)$-dimensional knots ($k\in\N\cup\{0\}$),  
we have an identity 
$$A(K_+)-A(K_-)=(t+1)A(K_0)$$ 
associated with the twist move.
It is a new type of local move identities of knot polynomial. 
Note $(t+1)$ in the right hand side. It is not $(t-1)$.

\bigbreak
\noindent
{\bf 
A combination of knot products and local moves on $n$-knots ($n\in\N$), and 
our main results.}  
In \cite{KauffmanOgasa, KauffmanOgasaII} 
we combined research of knot products and that of local moves on high dimensional knots,  
and obtained many results. 
We cite main results of them below.  
We generalize them and obtain our main new results, Theorems \ref{Bunkyo} and \ref{aoiro}, of this paper.\\  

On high dimensional pass-moves 
we prove the following Theorems \ref{Komagome}  and \ref{Bunkyo}. \\  

In \cite[Theorem 8.1]{KauffmanOgasa} 
we proved the following: 
If a 1-knot $A$ is obtained from a 1-knot $B$ by one pass-move, 
then 
$A\otimes^{\mu} {\rm Hopf}$ is obtained from $B\otimes^{\mu} {\rm Hopf}$
by one $(2\mu+1, 2\mu+1)$-pass-move.   \\

In \cite[Theorem 8.10]{KauffmanOgasa} 
we proved the case where $A$ is a knot in Theorem \ref{Komagome}.

\begin{thm}\label{Komagome} 
Let $A$ be a 1-link. 
Let $\mu\in\N$. 
Let $J=A\otimes^{\mu} {\rm Hopf}$. 
Let $K$ be obtained from $J$ by one $(2\mu+1, 2\mu+1)$-pass-move. 
Then there is a 1-link $B$ such that  $K=B\otimes^{\mu} {\rm Hopf}$  and such that 
 $A$ is pass-equivalent to $B$.  
\end{thm}

\cite[Theorems 8.1 and 8.10]{KauffmanOgasa} 
imply the following: 
A 1-knot $A$ is pass-equivalent to a 1-knot $B$ if and only if 
$A\otimes^{\mu} {\rm Hopf}$ is 
$(2\mu+1, 2\mu+1)$-pass-equivalent to 
$B\otimes^{\mu} {\rm Hopf}$.
We generalize this to the link case as follows.  
\\

In \cite[Main Theorem 4.2]{KauffmanOgasaII} 
we proved the following: 
If a 1-link $A$ is obtained from a 1-link $B$ by one pass-move, 
then 
$A\otimes^{\mu} {\rm Hopf}$ is obtained from $B\otimes^{\mu} {\rm Hopf}$
by one $(2\mu+1, 2\mu+1)$-pass-move  (see Note \ref{ei}: This fact is also proved by using Theorem \ref{Chicago}). 
By this fact and Theorem \ref{Komagome} we have the following.

\begin{thm}\label{Bunkyo}   Let  $\mu\in\N.$
A 1-link $A$ is pass-equivalent to a 1-link $B$ if and only if 
$A\otimes^{\mu} {\rm Hopf}$ is 
$(2\mu+1, 2\mu+1)$-pass-equivalent to 
$B\otimes^{\mu} {\rm Hopf}$.
\end{thm}

On twist-moves on high dimensional knots we prove the following Theorem \ref{aoiro}.  
\\

\cite[Theorem 4.1]{KauffmanOgasa} proved  the following: 
Suppose that two 1-links $J$ and $K$ differ by a single crossing change.
Then the knot products, 
$J \otimes^\mu{\rm Hopf}$ 
and   
$K\otimes^\mu{\rm Hopf}$,   
differ by a single twist-move, 
where $\mu\in\N\cup\{0\}$. 
See the left two figures of Figure \ref{Hawaii} for the crossing change on 1-links. 
\\

 \cite[Theorem 7.1]{KauffmanOgasa} proved the following: 
Let $m\in\N\cup\{0\}$. 
Suppose that two 
$($not necessarily connected$)$ $(2m+1)$-dimensional closed oriented 
submanifolds  of $S^{2m+3}$, 
$J$ and $K$,  differ by a single twist-move.
Then
the $(2m+2\nu+1)$-submanifolds of $S^{2m+2\nu+3}$, 
$J\otimes^\nu[2]$  
and 
$K\otimes^\nu[2]$, 
differ by a single twist-move.
Note that Hopf$=[2]\otimes[2]$
\\

 \cite[Theorem 7.3]{KauffmanOgasa} proved the following: 
Let $k\in\N$.    
Let $K$ $($respectively,  $J)$ be $(4k+5)$-submanifold of $S^{4k+7}$. 
Suppose that $K$ and $J$ differ by a single twist-move and are nonisotopic. 
Suppose that $K$ is isotopic to $A\otimes^{k+1}{\rm Hopf}$ for a 1-knot $A$. 
Then  there is a unique equivalence class of simple $(4k+1)$-knots for $K$ $($respectively,  $J)$   
with the following properties. 

\smallbreak
\noindent $\mathrm{(i)}$
There is a representative element $K'$ of the above equivalence class for $K$ such that 
$K$ is isotopic to $K'\otimes{\rm Hopf}$.   

\smallbreak
\noindent $\mathrm{(ii)}$
There is a representative element $J'$ of the above equivalence class for $J$ such that 
$J$ is isotopic to $J'\otimes{\rm Hopf}$.   

\smallbreak
\noindent $\mathrm{(iii)}$
$K'$ and $J'$ differ by a single twist-move and are nonisotopic. 
\\

Compare the above \cite[Theorem 7.3]{KauffmanOgasa} with 
the following Theorem \ref{aoiro}.

\begin{thm}\label{aoiro}     
Let $p\geq3$  
and $p\in\N$.    
Let $J$ be a $(2p+1)$-submanifold of $S^{2p+3}$. 
Suppose that $J$ and $K$ differ by a single twist-move and  are nonisotopic. 
Suppose that $J$ is isotopic to $A\otimes[2]$ for 
a $(2p-1)$-dimensional connected, $(p-2)$-connected, simple 
submanifold $A\subset S^{2p+1}$. 
Then  there is a $(2p-1)$-dimensional connected, $(p-2)$-connected, simple 
submanifold $B\subset S^{2p+1}$ with the following properties: 

\smallbreak
\noindent $\mathrm{(i)}$
$K$ is isotopic to $B\otimes[2]$.   

\smallbreak 
\noindent $\mathrm{(ii)}$
$A$ and $B$ differ by a single twist-move and are nonisotopic.

\smallbreak
\noindent $\mathrm{(iii)}$
The equivalence class of such $B$ is unique. 
\end{thm}

When we wrote about Theorem \ref{aoiro} in the abstract and \S\ref{Introduction}, 
we used the phrase `two-fold cyclic suspension'. 
It means taking a knot product with the empty knot $[2]$. 
\cite[Lemmas 5.3 and 5.4]{KauffmanNeumann} explain why we use it for the meaning.

In order to prove the above results in \cite{KauffmanOgasa, KauffmanOgasaII}, 
it is enough to use \cite{Levineun, Levinesimp}. 
However in order to generalize the above results in \cite{KauffmanOgasa, KauffmanOgasaII} 
and to obtain 
Theorems \ref{Komagome}, \ref{Bunkyo}, and \ref{aoiro}, 
we must generalize  \cite{Levineun, Levinesimp} and prove Theorem \ref{Chicago}.  
Theorem \ref{Chicago} implies Theorem \ref{carrot}, which is a result on Brieskorn submanifolds.

\bigbreak
\section{Review of simple submanifolds and that of Brieskorn submanifolds} \label{RkB}  
\noindent
In order to prove our main results, Theorems \ref{Bunkyo} and \ref{aoiro},  
we need to prove Theorem \ref{Chicago}.  
In order to prove Theorem \ref{Chicago} we review the following. 

\bigbreak

Let $K$ be a (not necessarily spherical) connected, closed, oriented, $(2p+1)$-submanifold of $S^{2p+3}$ $(p\in\N)$. 
$K$ is said to be {\it simple} if $K$ satisfies that \newline
$\pi_1(S^{2p+3}-N(K))\cong\Z$ and 
$\pi_i(S^{2p+3}-N(K))\cong0$ $(2\leqq i\leqq p)$. 
Let $V$ be a Seifert hypersurface for $K$. 
$V$ is said to be {\it simple} if $\pi_i(V)$ is trivial  $(1\leqq i\leqq p)$. 
See Theorem \ref{handle}. 
\\

Let $V$ be a Seifert hypersurface for a closed oriented $(2p+1)$-submanifold $K\subset S^{2p+3}$.    
Let $x_1,..., x_\mu$ be $(p+1)$-cycles in $V$ 
which compose an ordered basis of $H_{p+1}(V;\Z)$/Tor. 
Recall that the orientation of $V$ is compatible with that of $K$. 
Push $x_i$ into the positive direction of the normal bundle of $V$. 
Call it $x_i^{+}$.  
Push $x_i$ into the negative direction of the normal bundle of $V$. 
Call it $x_i^{-}$.  
A $(p+1)$-{\it $($positive$)$ Seifert matrix}   
for $K$ associated with $V$ represented by 
the ordered basis,  
$\{x_1,..., x_\mu\}$,  
is a $(\mu\x\mu)$-matrix 
$$S=(s_{ij})=({\mathrm{lk}}(x_i, x_j^{+})).$$

A $(p+1)$-{\it negative Seifert matrix} 
for $K$ associated with $V$ represented by 
the ordered basis,  
$\{x_1,..., x_\mu\}$,   
is a $(\mu\x\mu)$-matrix 
$$N=(n_{ij})=({\mathrm{lk}}(x_i, x_j^{-})).$$

We sometimes omit to write $K$, $V$, and $\{x_i\}$ 
when they are clear from the context. 
See e.g. \cite{KauffmanOgasa} for $(p,n+1-p)$-Seifert matrices for $n$-knots ($p, n\in\N$). 

\bigbreak
Let $A$ be an $r\x r$-matrix.  
Let $P$ be a unimodular $r\x r$-matrix. 
We say that $A$ is {\it equivalent} to $A'$ 
if $A'=^t\hskip-2mm PAP$, 
where $^t$ denotes an operation of making a transposed matrix.

\begin{pr}\label{tea} 
Let $S$ $($respectively,  $S')$ 
be a $(p+1)$-positive Seifert matrix 
for the above $(2p+1)$-knot $K$ associated with the above $V$
represented by an ordered basis $\{x_i\}$ $($respectively,  $\{x'_i\})$   
of $(p+1)$-cycles.  
Then $S$ is equivalent to $S'$.  
\end{pr}

\noindent{\bf Proof of Proposition \ref{tea}.}
Let $P$ be a matrix changing the ordered basis $\{x_i\}$ 
into the other $\{x'_i\}$. Then 
$S'=^t\hskip-2mm PSP$. 
(Someone may write 
$S'=PS (^t\hskip-1mm P)$,  
$S'=^t\hskip-2mm (P^{-1})S(P^{-1})$, or 
$S'=(P^{-1})S(^t (P^{-1}))$,  
in another convention.)
\qed\\

If a $(p+1)$-positive Seifert matrix $P$ for $K$ 
and a $(p+1)$-negative Seifert matrix $N$ for $K$ 
are defined by the same ordered basis $\{x_i\}$,    
we say that $P$ and $N$ are {\it related} and that the pair $(P, N)$ is a 
{\it pair of $(p+1)$-related Seifert matrices} for $K$.

\begin{pr}\label{daiji} 
Let $X$ be a $(p+1)$-positive Seifert matrix for a $(2p+1)$-dimensional closed oriented submanifold $K\subset S^{2p+3}$ associated with a Seifert hypersurface $V$. 
Then we have the following. 

\smallbreak\noindent$(1)$
$(-1)^p\hskip1mm^t\hskip-1mm X$ represents the $(p+1)$-negative Seifert matrix related to $X$ for $K$. 

\smallbreak\noindent$(2)$
$X-(-1)^p\hskip1mm^t\hskip-1mm X$ represents the intersection product 
on $H_{p+1}(V;\Z)${\rm /Tor}.  
\end{pr}

\noindent{\bf Proof of Proposition \ref{daiji}.}
By the definition of $y^{+}$, $x^{-}$, we have  
 ${\mathrm{lk}}(x, y^{-})={\mathrm{lk}}(x^{+}, y)$.    
By \cite[P 541]{LevineAlex} 
 ${\mathrm{lk}}(x^{+}, y)=(-1)^p{\mathrm{lk}}(y, x^{+})$.       
By these facts, 
the definition of $(p+1)$-positive Seifert matrices, and 
that of  $(p+1)$-negative Seifert matrices, 
we have Proposition \ref{daiji}.(1). 

Let $x$ and  $y$ be $(p+1)$-cycles. \\
${\mathrm{lk}}(x, y^{+})-(-1)^p{\mathrm{lk}}(y, x^{+})$ \\
$={\mathrm{lk}}(x, y^{+})-{\mathrm{lk}}(x^{+}, y)$\\
$={\mathrm{lk}}(x, y^{+})-{\mathrm{lk}}(x, y^{-})$\\
$={\mathrm{lk}}(x, y^{+}\amalg (-y^{-}))$, where $\amalg$ denotes the disjoint union. 
By using the $(p+2)$-chain $y\x[-1,1]$ which bounds $y^{+}\amalg (-y^{-})$,  
we can prove that 
${\mathrm{lk}}(x, y^{+}\amalg (-y^{-}))$ 
is the intersection product  $x\cdot y$ of $x$ and $y$ in $V$.  
%
%
Hence Proposition \ref{daiji}.(2) holds. 
\qed\\

Let $K$ be a closed oriented $(2p+1)$-submanifold of $S^{2p+3}$ $(p\in\N)$. 
If $A$ is a $(p+1)$-Seifert matrix (respectively,  negative $(p+1)$-Seifert matrix)
associated with a simple Seifert hypersurface  $V$ for $K$, 
we call  $A$ a {\it simple Seifert matrix}   
(respectively,  {\it negative simple Seifert matrix}
)
for $K$. 
We can define a 
{\it  pair of related simple Seifert matrices} 
for $K$. 
\smallbreak
\noindent
{\bf Note.} 
If we say that $K$ has a simple Seifert matrix $A$, 
then it means that $K$ has a simple Seifert hypersurface. 

\begin{defn}\label{KN} 
Let $*\in\{1,...,q\}.$ 
Let  $a_*$ be an integer$\geqq2$. 
Let $q\in\N$.  
The submanifold  \newline 
$\{(z_1,...,z_q)| |z_1|^2+,,,+|z_q|^2=1,  z_*\in\C\}\cap
\{(z_1,...,z_q)| z_1^{a_1}+...+z_q^{a_q}=0,  z_*\in\C\}$\newline 
$\subset\{(z_1,...,z_q)|  |z_1|^2+,,,+|z_q|^2=1, z_*\in\C\}$
is called the {\it Brieskorn submanifold $\Sigma(a_1,...,a_q)$.} 
The oriented diffeomorphism type of the Brieskorn submanifold $\Sigma(a_1,...,a_q)$ 
is called the {\it Brieskorn manifold $\Sigma(a_1,...,a_q)$.}

Let $a$ be an integer$\geqq2.$
Let 
$\Lambda_a=$
$(\zeta_{i,j})$ be an $(a-1)\x(a-1)$ matrix such that 

\smallbreak
\noindent
$(\zeta_{i,j})=$ 
$
\begin{cases}
1&\text{if $i=j$} \\
-1&\text{if $j=i+1$}  \\
0&\text{else,} 
\end{cases}
$ 
%
\quad\quad that is, \quad\quad 
%
$(\zeta_{i,j})=$ 
$
\begin{pmatrix}
1&-1&   &  &    \\
  & 1&-1&  &    \\
  &  &\cdot&\cdot   &    \\ 
  &  &   & 1&-1 \\   
  &  &   &   &1 
\end{pmatrix}
$.

\noindent 
Let 
$\Lambda_{a_1,...,a_q}=(-1)^\frac{(q-1)q}{2}\Lambda_{a_1}\otimes...\otimes\Lambda_{a_q}$. 
It is called 
a {\it Kauffman-Neumann-type}, or a {\it $KN$-type}. 
See \cite[the last few paragraphs of \S6]{KauffmanNeumann} 
for the definition of $KN$-types.

The Brieskorn submanifold  $\Sigma(a,b)\subset S^3$ is the torus $(a,b)$ knot (see \cite{Milnor}). 
We say that the Brieskorn submanifold  $\Sigma(a)\subset S^1$ is the empty knot of degree $a$ (see \cite{KauffmanNeumann}). 
\end{defn}

\begin{thm}\label{KNB} {\rm(\cite{Kauffman, KauffmanNeumann}.)} 
Let $q$ be an integer$\geqq3$. 
Let $*\in\{1,...,q\}$.   
Let $a_*$ be an integer\\$\geqq2$. 
Any Brieskorn submanifold $\Sigma(a_1,...,a_q)$ is 
a $($not necessarily spherical$)$ $(2q-3)$-dimensional connected, $(q-3)$-connected,  
simple submanifold of $S^{2q-1}$ 
such that $\Lambda_{a_1,...,a_q}$ is a simple 
Seifert matrix.    
\end{thm}

We prove that 
the converse of the $q\geqq4$ case of Theorem \ref{KNB} holds and 
that that of the $q=3$ case does not hold in general.  
See Theorems \ref{Chicago} and \ref{carrot}, and Note \ref{igai}.

\bigbreak
\section{Theorems on simple submanifolds and Brieskorn submanifolds}\label{TSB} 
\noindent
We explain Theorems \ref{Chicago} and \ref{carrot}.  
They are results on simple submanifolds and Brieskorn submanifolds. 
We use Theorem \ref{Tokyo} and prove them. 

\begin{thm}\label{Tokyo}
Let $p\in\N$.
Let $K$ be a closed oriented $(2p+1)$-dimensional  
connected, $(p-1)$-connected, simple submanifold of $S^{2p+3}$.  
Let $P$ be a simple 
 Seifert matrix 
for $K$. 
Then the following two conditions are equivalent. 

\smallbreak\noindent{\rm(i)}
$P'$ is a simple 
 Seifert matrix for $K$. 

\smallbreak\noindent{\rm(ii)}
$P'$ is $(-1)^p$-$S$-equivalent to $P$.  
\end{thm}

Recall that the term `$(-1)^p$-$S$-equivalent' is defined in Definition \ref{saisho}.

Note that we have the following proposition. 
\begin{pr}\label{Osaka}
There is a natural number $p$ and a closed oriented $(2p+1)$-dimensional   
connected, $(p-1)$-connected, simple submanifold $K$ of $S^{2p+3}$
with the following property:  
There is a simple Seifert matrix $P$ and a Seifert matrix $R$ for $K$ such that 
$P$ is not $(-1)^p$-$S$-equivalent to $R$.  
\end{pr}

\noindent
{\bf Note.}
In Proposition \ref{Osaka},  $R$ is not a simple 
Seifert matrix. 
If $R$ is associated with a Seifert hypersurface $W$ for $K$, 
then $W$ is not a simple Seifert hypersurface.

\bigbreak 
It was proved in \cite{Levinesimp} 
that if $K$ and $J$ are spherical, 
Theorem \ref{Chicago} is true.   
We generalize the result in \cite{Levinesimp}
and prove a stronger theorem, which is Theorem \ref{Chicago}.
See the following notes to Theorem \ref{Chicago}.


\begin{note}\label{tsuika}   
\cite[the results in the first page]{Durfee} and 
\cite[the results in the first page]{MK} 
claim the following: 
Add the condition `$K$ and $J$ are fibered knots' to Theorem \ref{Chicago}.(1). 
Then Theorems \ref{Chicago}.(1).(i) and (ii) are equivalent under this condition. 

On the other hand, in Theorem \ref{Chicago}.(1), 
we do not suppose that $K$ is fibered nor that $J$ is fibered. 
Therefore our result, Theorem \ref{Chicago}.(1), is stronger than the results in \cite{Durfee, MK}.  
\end{note}

\begin{note}\label{medaka}   
\cite[Theorem 3.1]{Sa} claims the following: 
There is a pair $(K, J)$ with the following property: 
$K$ and $J$ are closed oriented 3-dimensional simple submanifolds of $S^5$. 
$K$ is diffeomorphic to $J$.   
$K$ is nonisotopic to $J$.   
There is a simple Seifert matrix $P_K$ (respectively,  $P_J$) for $K$ (respectively,  $J$)  
such that $P_K$ is $(-1)$-$S$-equivalent to $P_J$. 

Therefore we have 
neither the $p=1$ version of Theorem \ref{Chicago}.(1) nor that of Theorem \ref{carrot}.(1) in general. 
\end{note}

\begin{note}\label{tagame}
\cite[Theorem 2.2]{Sa} claims the following: 
Add the condition `$K$ and $J$ are fibered knots' to Theorem \ref{Chicago}.(2). 
Then Theorems \ref{Chicago}.(2).(i) and (ii) are equivalent under this condition. 

On the other hand, in Theorem \ref{Chicago}.(2), 
we do not suppose that $K$ is fibered nor that $J$ is fibered. 
Therefore our result, Theorem \ref{Chicago}.(2), is stronger than 
\cite[Theorem 2.2]{Sa}.  

\cite[Theorem 2.2]{Sa0} claims the following without a proof: 
Add the condition `$K$ and $J$ are diffeomorphic to a homology 3-sphere' 
to Theorem \ref{Chicago}.(2).  
Then Theorems \ref{Chicago}.(2).(i) and (ii) are equivalent under this condition. 

On the other hand, in Theorem \ref{Chicago}.(2), 
we do not suppose 
that $K$ is diffeomorphic to a homology 3-sphere 
nor
that $J$ is diffeomorphic to a homology 3-sphere. 
Therefore our result, Theorem \ref{Chicago}.(2), is stronger than 
\cite[Theorem 2.2]{Sa0}.  
\end{note}

\noindent{\bf Note.}   
 Theorem \ref{Chicago}.(2) is not trivial. Reason: 
If we assume a similar condition for 1-knots (respectively,  1-links)  
to Theorem \ref{Chicago}.(2).(ii), that is, 
replace 2-handles (respectively,  2-cycles) with 1-handles (respectively,  1-cycles),  
then the 1-knots (respectively,  1-links), $K$ and $J$, which we obtain, are nonisotopic in general.  

By using \cite{Wall1, Wall2}, we could 
sophisticate Theorem \ref{Chicago}.(2).(ii).

\begin{note} \label{igai}
On Theorem \ref{Chicago}.(2) and the paragraph under Theorem \ref{KNB}, note the following. 
There is a closed oriented simple 3-submanifold $E\subset S^5$ with the following properties: 

\smallbreak\noindent(i)   There is a simple Seifert matrix $P$ for $E$ which is $(-1)$-$S$-equivalent to the empty matrix. 

\smallbreak\noindent(ii)
$E$ is not homeomorphic to the 3-sphere.

\smallbreak\noindent See \cite[Note to Proof of Theorem 7.2.(1)]{KauffmanOgasa}. 
\end{note}
 
Theorem \ref{Chicago} and Definition \ref{KN} imply that 
the converse of the $q\geqq4$ case of Theorem \ref{KNB} holds. 
See 
Theorem \ref{carrot}.(1). Note $q=p+2$.  
Of course Theorem \ref{Chicago}.(2) holds if $J$ is the Brieskorn submanifold 
which is a 3-manifold. 
It is  Theorem \ref{carrot}.(2).  
See the following note to Theorem \ref{carrot}.  



\begin{note}\label{tsuikaBr}
On Theorem \ref{carrot}.(1), we have a similar situation to Note \ref{tsuika} as follows.   
The results in \cite[the first page]{Durfee} and 
the results in \cite[the first page]{MK} 
claim the following: 
Add the condition `$K$ and $J$ are fibered knots' to Theorem \ref{carrot}.(1). 
Then Theorems \ref{carrot}.(1).(i) and (ii) are equivalent under this condition. 

On the other hand, in Theorem \ref{carrot}.(1), 
we do not suppose that $K$ is fibered nor that $J$ is fibered. 
Therefore our result, Theorem \ref{carrot}.(1), is stronger than the results in \cite{Durfee, MK}.  
\end{note}


By Theorem \ref{Chicago} we have the following.

\begin{thm}\label{Boston} 
$(1)$ 
Let $p\geq2$ 
and $p\in\N$. 
Let $K$ be a closed oriented $(2p+1)$-submanifold of $S^{2p+3}$. 
Let $A$ be a 1-link. 
Let $a_*$ be an integer $\geq2$ $(*=1,...,p)$. 
Then the following two statements are equivalent.

\smallbreak\noindent{\rm(i)}
$K=A\otimes[a_1]\otimes...\otimes[a_p]$. 

\smallbreak\noindent{\rm(ii)}
$K$ is connected, $(p-1)$-connected, and simple. 
There is a simple 
Seifert matrix $P$ for $K$  
and a Seifert matrix $P'$ for $A$ 
such that $P$ is $(-1)^p$-$S$-equivalent to 
$P'\otimes\Lambda_{a_1,...,a_p}$.

\bigbreak\noindent$(2)$     
Let $p\geq2$ 
and $p\in\N$. 
Let $K$ be a closed oriented $(2p+1)$-submanifold of $S^{2p+3}$. 
Let $p>q.$ 
Let $2q+1\geq3$ and $q\in\N$. 
Let $A$ be a simple $(2q+1)$-submanifold.  
Let $a_*$ be an integer $\geq2$ $(*=1,...,p-q)$. 
Then the following two statements are equivalent.

\smallbreak\noindent{\rm(i)}
$K=A\otimes[a_1]\otimes...\otimes[a_{p-q}]$. 

\smallbreak\noindent{\rm(ii)}
$K$ is connected, $(p-1)$-connected, and simple.   
There is a simple 
Seifert matrix  $P$ for $K$ 
and a simple 
Seifert matrix $P'$ for $A$ such that 
$P$ is $(-1)^p$-$S$-equivalent to \newline
$(-1)^{q(p-q)}P'\otimes\Lambda_{a_1,...,a_{p-q}}$.  
\end{thm}

\begin{note}\label{MA}  
In 
Theorem \ref{Boston} 
 we use the following fact, which is proved in \cite[\S6]{KauffmanNeumann}: 
Let $X$ (respectively,  $Y$) be a $(2x+1)$-(respectively,  $(2y+1)$-)dimensional  
closed oriented simple submanifold of $S^{2x+3}$ (respectively,  $S^{2y+3}$), 
where $x,y\in\N\cup\{-1,0\}$. Here,  
we regard $X$ (respectively,  $Y$) as a 1-link if $x$ (respectively,  $y$) is 0,   
and 
as the empty knot if $x$ (respectively,  $y$) is $-1$.  
Let $S_X$ (respectively,  $S_Y$) be 
an $(x+1)$-(respectively,  $(y+1)$-)positive Seifert matrix for $X$ (respectively,  $Y$). 
Let \newline
\hskip3cm$S_{X\otimes Y}=(-1)^{xy}S_X\otimes S_Y$.   \newline
Then $S_{X\otimes Y}$ is an $(x+y+3)$-positive Seifert matrix for $X\otimes Y$.  
Furthermore it holds that $X\otimes Y$ is a $(2x+2y+5)$-dimensional  
closed oriented simple submanifold of $S^{2x+2y+7}$ 
if 
$x+y\geqq-1$. 
\end{note}
\bigbreak

Let $M$ be a closed oriented $m$-manifold which we can  embed in $S^{m+2}$.   
Let $K$ be a closed oriented $m$-submanifold of $S^{m+2}$ which is diffeomorphic to $M$. 
Let $K$ be an image of a smooth embedding map $f:M\to S^{m+2}$. 
Take a tubular neighborhood $N(K)$ of $K$ in $S^{m+2}$. 
Note that it is diffeomorphic to $M\x D^2$. Of course 
$\partial(M\x D^2)=M\x \partial D^2\\=M\x S^1$. 
Let $N(K)$ be an image of a smooth embedding map $F:M\x D^2\to S^{m+2}$.  

Let $h:M\x S^1\to M\x S^1$ be the diffeommorphism defined by 
$(x,y)\mapsto (\Psi(x)\cdot x, y)$ where $\Psi:S^1\to$Diff$_+$ $M$. 
Here,  Diff$_+$ $M$ is a set of orientation preserving diffeomorphism maps of $M$.

\begin{que}\label{qichi}
Attach $M\x D^2$ with 
$\overline{S^{m+2}-N(K)}$ by a diffeomorphism $F\circ h$. 
Suppose that 
$(M\x D^2)\cup_{F\circ h}\overline{S^{m+2}-N(K)}$ is diffeomorphic to $S^{m+2}$. 
Thus we obtain a submanifold $M\x\{0\}\subset S^{m+2}$, say $K_f$. 
Then is $K_f$ isotopic to $K$?  
\end{que}

\cite[section 22]{Levinesimp} gave a partial solution to Question \ref{qichi}. 
It implies the following: 
Let $m=2p+1$,   $p\geq2$ and $p\in\N$. 
Let $K$ be a spherical simple knot. 
Then we have an affirmative answer to Question \ref{qichi}. 

Its generalization is obtained by Theorem \ref{Chicago}.(1). 
It is also a partial solution to Question \ref{qichi}.
That is, we have the following: 

\begin{cor}\label{hokuukan}  
Let $p\geq2$ 
and $p\in\N$.  
In Question \ref{qichi}, if $K$ is a closed oriented $(2p+1)$-dimensional 
connected, $(p-1)$-connected, simple submanifold of $S^{2p+3}$, we have the positive answer. 
\end{cor}

Of course, Question \ref{qichi} is motivated by the following question.

\begin{que}\label{qni}
In Question \ref{qichi},  replace $F\circ h$ with any orientation preserving diffeomorphism of 
$M\x S^1$. 
What is the answer? 
 
That is, in other words, 
is the submanifold type of $K$ determined by the diffeomorphism type of the complement
$\overline{S^{m+2}-N(K)}$?
\end{que}

A combination of \cite[section 22]{Levinesimp} and \cite{Gluck} implied the following:
Let $K$ be a simple spherical $(2p+1)$-knot ($p\geq2$). 
Let $J$ be a spherical   $(2p+1)$-knot. 
If $\overline{S^{2p+3}-N(K)}$ and $\overline{S^{2p+3}-N(J)}$ are diffeomorphic, 
$K$ and $J$ are isotopic.

\cite{CS} proved that, 
if $m=3,4$, 
spherical knots are not determined by the complements.

\cite[Theorem 4.6]{Sa} proved the following. 
There are simple 3-dimensional submanifolds $K$ and $J$ of $S^5$ with the following properties: 
$K$ and $J$ are diffeomorphic as abstract manifolds.  
A simple $S$-matrix of $K$ and that of $J$ are $(-1)$-$S$-equivalent. 
$\overline{S^5-N(K)}$ and $\overline{S^5-N(J)}$  are diffeomorphic as abstract manifolds. 
$J$ and $K$ are nonisotoopic.

There are many other important results on Questions \ref{qichi} and \ref{qni}.  
See \cite{Gordon, GL, Suciu} etc. 
We have not answered Questions \ref{qichi} nor \ref{qni} completely in all cases.

\bigbreak
\section{Handle decompositions of Seifert hypersurfaces} \label{handleproof}
\noindent
We use Theorem \ref{handle} in order to prove our other theorems in this paper. 
We need preliminaries before stating it.
\\

We review a few facts on handle decompositions. 
See \cite{Browder, Kirby, Luck, Ranichi, Smale, Wall} for detail.

Let $W$ be a $w$-dimensional compact manifold ($w\in\N$). 
Take a handle decomposition \newline  
$(B\x[0,1])\cup$(0-handles)$\cup$...$\cup$($w$-handles)$\cup(T\x[0,1])$, 
where there may not be an $i$-handle ($0\leqq i\leqq w$). 
Note that 
$B$ (respectively,  $T$) is a compact $(w-1)$-submanifold of $\partial W$. 
Note that  
$\partial B$ is diffeomorphic to $\partial T$.  
Note that 
$\partial W$ is diffeomorphic to 
$B\cup_\alpha T$, where $\alpha$ is a diffeomorphism map from 
$\partial B\to\partial T$.   
%
%
If a handle is attached to $B\x[0,1]$, its attaching part is embedded in $B\x\{1\}$. 
No handle is attached to $T\x[0,1]$ although \\(a handle)$\cap(T\x[0,1])\neq\phi$ may hold.
If (a handle)$\cap(T\x[0,1])\neq\phi$, the intersection\\$\subset T\x\{0\}.$ 
$B$ (respectively,  $T$) may be the empty set. 
We do not suppose whether $B$ (respectively,  $T$) is connected or not.  
We do not suppose whether $B$ (respectively,  $T$) is closed or not.  
We say that $B$ or $B\x\{0\}$ is the {\it bottom} of this handle decomposition 
and that $T$ or $T\x\{1\}$ is the {\it top} of this handle decomposition. \\


For a handle decomposition, 
 regard the core (respectively,  cocore) of each handle as the cocore (respectively,  core) of it 
and replace the top (respectively,  the bottom) with the bottom (respectively,  the top), 
then we obtain a new  handle decomposition. It is called 
the {\it dual handle decomposition} of the handle decomposition. 
If a handle $h$ in a handle decomposition is changed into  
a handle $\bar{h}$ in its dual handle decomposition, 
$\bar{h}$ is called the {\it dual handle} of $h$.  
Note that we have the following: 
The top of this dual handle decomposition is $B\x\{0\}$, 
and the bottom of it is $T\x\{1\}$. 
If a dual handle is attached to $T\x[0,1]$, its attaching part is embedded in $T\x\{0\}$.  
No dual handle is attached to $B\x[0,1]$ 
although (a handle)$\cap(B\x[0,1])\neq\phi$ may hold.
If (a handle)$\cap(B\x[0,1])\neq\phi$, the intersection$\subset B\x\{1\}$.
\\

Let $V$ be a compact $(n+1)$-manifold ($n+1\in\N$). 
For the convenience of the application (see Theorem \ref{handle}),   
we suppose that the dimension is $n+1$, not $n$.  
If a handle decomposition of $V$ satisfies the following conditions, 
we say that the handle decomposition is a {\it special handle decomposition} of $V$. 

\smallbreak\noindent(1) 
The top $T$ is connected or empty. 
The bottom $B$ is connected or empty.

\smallbreak\noindent(2) 
It has only one (respectively,  no) $(n+1)$-dimensional 0-handle if $B=\phi$ (respectively,  $B\neq\phi$). 
It has no $(n+1)$-dimensional $i$-handle $(1\leqq i\leqq[\frac{n-1}{2}]).$  

\smallbreak\noindent(3)   
The dual handle decomposition has 
only one (respectively,  no) $(n+1)$-dimensional 0-handle if $T=\phi$ (respectively,  $T\neq\phi$). 
The dual handle decomposition has no  $(n+1)$-dimensional $i$-handle 
$(1\leqq i\leqq[\frac{n-1}{2}]).$  
\\

\noindent{\bf Example.}  
If the above $V$ is 6-dimensional, and 
has a special decomposition with $B=\phi$ and $T=\partial V$, 
then $V$ has a handle decomposition 
(one 0-handle)$\cup$(3-handles), where 
there may be no 3-handle. 
If  the above $V$ is 7-dimensional and 
has a special decomposition $B=\phi$ and $T=\partial V$, 
then $V$ has a handle decomposition  \newline 
(one 0-handle)$\cup$(3-handles)$\cup$(4-handles), 
where there may be no 3-handle, or where there may be no 4-handle. 
\\

We use 
`surgeries by using embedded handles' defined in Definition \ref{subsur}. 


\bigbreak
Let $n\in\N\cup\{0\}$. 
Let $K$ be an $n$-dimensional oriented closed submanifold of 
a (not necessarily closed) $(n+2)$-dimensional oriented compact manifold-with-boundary $Q$. 
We suppose that $K$ satisfies the following condition 
$(\star)$. 

\begin{figure}
\unitlength 0.1in
\begin{picture}(38.35,20.25)(1.30,-21.95)
%
\special{pn 8}%
\special{pa 1390 1670}%
\special{pa 1380 1640}%
\special{pa 1375 1609}%
\special{pa 1376 1577}%
\special{pa 1379 1545}%
\special{pa 1385 1513}%
\special{pa 1391 1480}%
\special{pa 1396 1448}%
\special{pa 1398 1416}%
\special{pa 1399 1384}%
\special{pa 1400 1352}%
\special{pa 1401 1320}%
\special{pa 1404 1288}%
\special{pa 1408 1257}%
\special{pa 1414 1225}%
\special{pa 1421 1194}%
\special{pa 1430 1164}%
\special{pa 1441 1133}%
\special{pa 1451 1103}%
\special{pa 1462 1073}%
\special{pa 1474 1043}%
\special{pa 1484 1012}%
\special{pa 1494 982}%
\special{pa 1506 952}%
\special{pa 1522 925}%
\special{pa 1541 899}%
\special{pa 1562 874}%
\special{pa 1584 851}%
\special{pa 1607 829}%
\special{pa 1632 808}%
\special{pa 1659 791}%
\special{pa 1687 776}%
\special{pa 1717 765}%
\special{pa 1748 757}%
\special{pa 1780 751}%
\special{pa 1812 747}%
\special{pa 1844 744}%
\special{pa 1876 742}%
\special{pa 1908 741}%
\special{pa 1939 741}%
\special{pa 1971 740}%
\special{pa 2003 740}%
\special{pa 2034 740}%
\special{pa 2067 741}%
\special{pa 2099 740}%
\special{pa 2132 740}%
\special{pa 2165 741}%
\special{pa 2194 750}%
\special{pa 2220 770}%
\special{pa 2243 793}%
\special{pa 2265 817}%
\special{pa 2286 841}%
\special{pa 2303 867}%
\special{pa 2318 895}%
\special{pa 2330 924}%
\special{pa 2340 954}%
\special{pa 2347 985}%
\special{pa 2354 1017}%
\special{pa 2360 1050}%
\special{pa 2365 1082}%
\special{pa 2371 1115}%
\special{pa 2377 1147}%
\special{pa 2382 1179}%
\special{pa 2388 1211}%
\special{pa 2393 1242}%
\special{pa 2399 1274}%
\special{pa 2404 1306}%
\special{pa 2408 1337}%
\special{pa 2412 1368}%
\special{pa 2416 1400}%
\special{pa 2419 1431}%
\special{pa 2421 1462}%
\special{pa 2422 1493}%
\special{pa 2423 1525}%
\special{pa 2425 1558}%
\special{pa 2426 1592}%
\special{pa 2425 1622}%
\special{pa 2421 1650}%
\special{pa 2409 1665}%
\special{pa 2374 1667}%
\special{pa 2346 1655}%
\special{pa 2317 1647}%
\special{pa 2286 1643}%
\special{pa 2254 1644}%
\special{pa 2221 1646}%
\special{pa 2188 1650}%
\special{pa 2154 1654}%
\special{pa 2121 1657}%
\special{pa 2089 1658}%
\special{pa 2057 1658}%
\special{pa 2025 1658}%
\special{pa 1993 1658}%
\special{pa 1962 1658}%
\special{pa 1930 1660}%
\special{pa 1898 1662}%
\special{pa 1866 1665}%
\special{pa 1834 1667}%
\special{pa 1802 1666}%
\special{pa 1771 1661}%
\special{pa 1739 1654}%
\special{pa 1708 1646}%
\special{pa 1677 1640}%
\special{pa 1645 1637}%
\special{pa 1613 1637}%
\special{pa 1581 1639}%
\special{pa 1549 1641}%
\special{pa 1517 1644}%
\special{pa 1485 1646}%
\special{pa 1453 1648}%
\special{pa 1422 1656}%
\special{pa 1391 1665}%
\special{pa 1390 1665}%
\special{sp}%
%
\special{pn 8}%
\special{pa 3705 1300}%
\special{pa 3703 1266}%
\special{pa 3720 1240}%
\special{pa 3737 1214}%
\special{pa 3748 1183}%
\special{pa 3758 1152}%
\special{pa 3771 1122}%
\special{pa 3784 1093}%
\special{pa 3797 1063}%
\special{pa 3808 1034}%
\special{pa 3816 1004}%
\special{pa 3821 972}%
\special{pa 3823 940}%
\special{pa 3824 908}%
\special{pa 3825 876}%
\special{pa 3826 844}%
\special{pa 3829 812}%
\special{pa 3831 780}%
\special{pa 3833 749}%
\special{pa 3832 717}%
\special{pa 3829 684}%
\special{pa 3823 651}%
\special{pa 3833 621}%
\special{pa 3873 600}%
\special{pa 3891 590}%
\special{pa 3857 591}%
\special{pa 3823 589}%
\special{pa 3805 569}%
\special{pa 3791 539}%
\special{pa 3773 512}%
\special{pa 3758 484}%
\special{pa 3746 453}%
\special{pa 3736 421}%
\special{pa 3728 387}%
\special{pa 3720 354}%
\special{pa 3712 321}%
\special{pa 3702 290}%
\special{pa 3690 261}%
\special{pa 3675 235}%
\special{pa 3655 214}%
\special{pa 3630 197}%
\special{pa 3600 185}%
\special{pa 3568 177}%
\special{pa 3535 172}%
\special{pa 3502 170}%
\special{pa 3490 170}%
\special{sp}%
\put(39.6500,-6.3500){\makebox(0,0)[lb]{$Q$}}%
\put(34.6000,-17.2000){\makebox(0,0)[lb]{$\partial Q$}}%
%
\special{pn 8}%
\special{pa 1135 1315}%
\special{pa 130 2195}%
\special{fp}%
%
\special{pn 8}%
\special{pa 130 2190}%
\special{pa 2785 2195}%
\special{fp}%
%
\special{pn 8}%
\special{pa 2790 2195}%
\special{pa 3695 1310}%
\special{fp}%
%
\special{pn 8}%
\special{pa 1140 1315}%
\special{pa 1335 1310}%
\special{fp}%
%
\special{pn 8}%
\special{pa 1465 1315}%
\special{pa 2330 1305}%
\special{fp}%
%
\special{pn 8}%
\special{pa 2475 1315}%
\special{pa 3695 1310}%
\special{fp}%
\put(12.7500,-8.3000){\makebox(0,0)[lb]{$K$}}%
\end{picture}%
\caption{{\bf An example of $K$ and $Q$.}\label{Delaware}}
\end{figure}
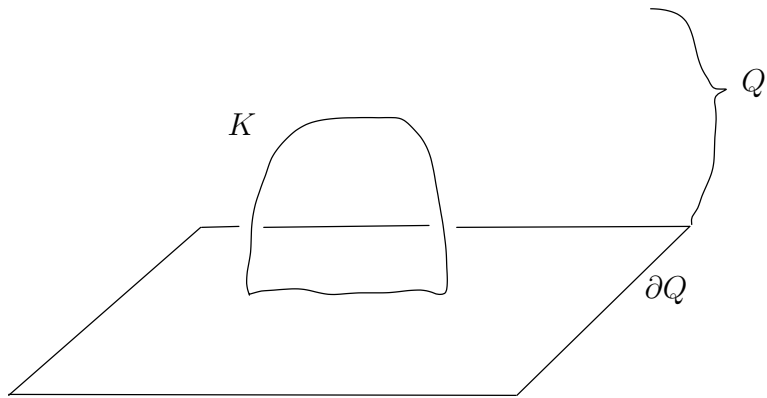

\bigbreak\noindent
{\bf The condition $(\star)$.}
$K\cap$ Int$Q=K-\partial Q$ is connected,   
and is an $n$-dimensional open submanifold of $K$.   
$\overline{K-\partial Q}$ is transverse to $\partial Q$, and 
is an $n$-dimensional compact submanifold of $Q$.  
$K\cap\partial Q$ is 
a (not necessarily connected) $n$-dimensional compact submanifold of $\partial Q$.
See Figure \ref{Delaware}  
for an example. 

\bigbreak
If $K\cap\partial Q\neq\phi$, we define the tubular neighborhood $N(K)$ of $K$ in $Q$ as follows. 
Let $N(K)$ be diffeomorphic to $K\times D^2$. 
Take the tubular neighborhood of  $\overline{K-\partial Q}$ in $Q$, and say $X$. 
Suppose that $\overline{N(K)-X}$ is as follows. 
Take the tubular neighborhood of 
$(K\cap\partial Q)-{\rm Int}X$ in $\partial Q$, and say $Y$. 
Take the unit inward vector at each point of $Y$. 
Note that the union of these vectors is $Y\x[0,1]$ in $Q$, and say $Z$.   
Let $\overline{N(K)-X}$ be $Z$. That is,  $N(K)=X\cup Z$. 



\begin{thm}\label{handle}   
Let $n\in\N$ and $n\geqq3$.   
Let $K$ be a closed oriented $n$-dimensional connected, $([\frac{n}{2}]-1)$-connected submanifold of $Q$.  
Let $Q$ be $S^{n+2}$,  $B^{n+2}$, or $S^{n+1}\x[0,1]$.   
We suppose that $K$ satisfies the above condition $(\star)$. 
We do not suppose whether $K\cap \partial Q=\phi$ or $K\cap \partial Q\neq\phi$.    
Let $N(K)$ be the tubular neighborhood of $K$ in $Q$. 
Then the following three conditions are equivalent: 

\smallbreak
\noindent {\rm (i)}   $\pi_i(\overline{Q-N(K)})=
\begin{cases}
\Z & \text {if $i=1$} \\
0  & \text{ if $2\leqq i\leqq[\frac{n-1}{2}]$.}
\end{cases}
$

\smallbreak
\noindent {\rm (ii)}
There is a $[\frac{n-1}{2}]$-connected Seifert hypersurface $V$ for $K$.  
$($Note: $V$ is $(n+1)$-dimensional.$)$

\smallbreak
\noindent {\rm (iii)} 
There is a Seifert hypersurface $V$ for $K$ that has a special handle decomposition 
whose bottom is the empty set and whose top is $\partial V$. 
\end{thm}  

\noindent{\bf Note.}
\cite{Levineun, Levinesimp} proved 
a special case of Theorem \ref{handle},   
which is the case where $K$ is PL homeomorphic to the standard sphere 
and $Q$ is $S^{n+2}$.  
Theorem \ref{handle} is its generalization, and is stronger than this special case. 
\\

\noindent{\bf Proof of Theorem \ref{handle}.}
Since $[K]=0\in H_n(Q;\Z)$, 
there is a Seifert hypersurface $V$ for $K$ 
even if $K\cap\partial Q\neq\phi$.  
By using isotopy of $V$, we can suppose that Int $V$ $\subset$ Int $Q$ 
and that $V$ is transverse to $\partial Q$. 
Let $N(V)$ be the tubular neighborhood of $V$ in $Q$.   
Note that 
$\overline{\overline{Q-N(K)}-N(V)}$ is isotopic to 
$\overline{Q-N(V)}$ in $Q$.

We prove the following proposition. 
\begin{pr}\label{ca}
`$({\rm iii})\Longrightarrow({\rm i})$' in Theorem $\ref{handle}.$  
\end{pr}

\noindent{\bf Proof of Proposition \ref{ca}.} 
By van Kampen theorem $\pi_1(\overline{Q-N(V)})\cong1$. 
Reason: Note that $Q=N(V)\cup\overline{Q-N(V)}$,  
and that $N(V)\cap\overline{Q-N(V)}$ is a disjoint union of two copies of $V$.  
Take a 1-handle $h^1$ embedded in $\overline{Q-N(V)}$, and attach it to $N(V)$
so that the two attached parts are in different components. 
Use van Kampen theorem to  a pair ($N(V)\cup h^1$, $\overline{Q-N(V)}$). 

By Mayer-Vietoris exact sequence on  
$N(V)$ and $\overline{Q-N(V)}$, 
we have 
$H_i(\overline{Q-N(V)};\Z)\\\cong0$ for $2\leq i\leq[\frac{n-1}{2}]$. 

Consider the infinite cyclic covering space $\widetilde{\overline{Q-N(K)}}$
of $\overline{Q-N(K)}$. 
Note that  $\widetilde{\overline{Q-N(K)}}$ is 
a union of the lift of $N(V)$ and that of $\overline{Q-N(V)}$.  
By Mayor-Vietoris exact sequence on them, 
$H_i(\widetilde{\overline{Q-N(K)}};\Z)\cong0$ if $1\leq i\leq[\frac{n-1}{2}]$.
By van Kampen theorem $\pi_1(\widetilde{\overline{Q-N(K)}})\cong1$. 
By Hurewicz theorem  
$\pi_i(\widetilde{\overline{Q-N(V)}})\cong0$ for $2\leq i\leq[\frac{n-1}{2}]$. 

By well-known facts on covering spaces, 
$\pi_i(\widetilde{\overline{Q-N(V)}})\cong\pi_i(\overline{Q-N(V)})$ 
for \newline$2\leq i\leq[\frac{n-1}{2}]$. 
This completes the proof of Proposition \ref{ca}. \qed

\begin{pr}\label{ab}
`$({\rm i})\Longrightarrow({\rm ii})$' in Theorem $\ref{handle}.$    
\end{pr}

\noindent{\bf Proof of Proposition \ref{ab}.}
 We prove the following claim.

\begin{cla}\label{mazu}
There is a simply-connected Seifert hypersurface for $K$. 
\end{cla}

\noindent{\bf Proof of Claim \ref{mazu}.} 
Take a Seifert hypersurface $V$ for $K$. 
Let $\{g_1,...,g_\mu\}$ be a set of generators of $\pi_1V$. 
Let $g_i$ also denote a circle that represents the element $g_i\in\pi_1V.$
Note that the dimension of $V$ is $n+1$ and that $n+1\geqq4$. 
Hence we can suppose that $g_i$ is embedded in $V$. 

Since $g_i$ is embedded in $V$, 
the intersection product of 
$[g_i]\in H_1(\overline{Q-N(K)}; \Z)$ 
and 
$[V,K]\in H_{n+1}(\overline{Q-N(K)}; \partial N(K);\Z)$  
is zero. 
Furthermore recall $\pi_1(\overline{Q-N(K)})=\Z$. 
Hence we can take a continuous map $f_i:D^2_i\rightarrow \overline{Q-N(K)}$ 
such that $f_i(\partial D^2_i)=g_i$.  

Note that the dimension of $\overline{Q-N(K)}$ is  $n+2$, 
and that $n+2\geqq5$. Hence we can suppose that $f_i$ is an embedding map. 

We need the following lemma.  

\begin{lem}\label{circle} 
Let $C$ be an embedded circle in $V$. 
Suppose that there is an embedded 2-disc $D$ in $\overline{Q-N(K)}$ 
such that $\partial D=C$. 
Note that $D$ may intersect $V$. 
Let $N(C)$ be the tubular neighborhood of $C$ in $\overline{Q-N(K)}$. 
Then we can suppose that $D\cap V\cap N(C)$ is $C$. 
$($Recall that $n+2\geqq5.)$
\end{lem}

\noindent{\bf Proof of Lemma \ref{circle}.}  
Since $V$ is orientable,  $N(C)=C\x D^{n+1}$.  
Hence $\partial N(C)$ is the trivial $S^n$-bundle over $C$. 
Since $n\geq2$, all sections of this trivial $S^n$-bundle over $C$ are homotopic. 
There is a section perpendicular to $V$ at $C$. 
Another section is defined by $D\cap N(C)$. 
Both sections are homotopic.  
Hence   $D\cap V\cap N(C)=C$. 
This completes the proof of Lemma \ref{circle}. \qed \bigbreak

Let $N(g_i)$ be the tubular neighborhood of $g_i$ in $\overline{Q-N(K)}$. 
By Lemma \ref{circle} \newline
$f_i(D^2_i)\cap V\cap N(g_i)=g_i$. 
Hence we can suppose 
that $f_i(D^2_i)$ intersects $V$ transversely,    
that $(f_i(D^2_i)\cap V)-g_i$ is a disjoint union of some circles,  
and that $f^{-1}$(the circles) is in the interior of $D^2_i$. 
Take an innermost circle of the circles $f_i(D^2_i)\cap V\subset D^2_i$ .   
It bounds a disc in $D^2_i$. Note that (the disc) $\cap$ (the other circles)$=\phi$. 
Note that $f($the disc$)$ is embedded in $\overline{Q-N(K)}$ and that  
 $f($the disc$)\cap V=f($the innermost circle$)$.  
Take an embedded $(n+2)$-dimensional 2-handle whose core is $f($the disc$)$ and 
whose attaching part is embedded in $V$.   
Carry out a surgery on $V$ by using this 2-handle.  
(Note that we remove the interior of $g_i\x D^n$ and add $D^2\x S^{n-1}$.)  
Repeating this procedure on such circles in each $D_i$, we obtain  a new $V$. 
By van Kampen theorem this new $V$ is simply-connected. 
This completes the proof of Claim \ref{mazu}. \qed

\bigbreak
It is trivial that Proposition \ref{ab} follows from the following Claims \ref{pi} and \ref{piinj}. 

\begin{cla}\label{pi}
Let $N(V)$ be the tubular neighborhood of  $V$ in $\overline{Q-N(K)}$. 
Note that $N(V)\newline=V\x [-1,1].$  
Let $r\leqq[\frac{n-1}{2}]$. 
Suppose that there is an $(r-1)$-connected Seifert hypersurface. 
Then there is an $(r-1)$-connected Seifert hypersurface $V$ with the following condition: 
For $t=1,-1$,  the homomorphism 
$\iota:\pi_r(V\x\{t\})\to\pi_r(\overline{\overline{Q-N(K)}-N(V)})$ 
that is induced by the natural inclusion map is injective.  
\end{cla}  


\noindent{\bf Proof of Claim \ref{pi}.}  
The $r=1$ case follows from Claim \ref{mazu}.  

We prove the $r\geqq2$ case.    
Take an $(r-1)$-connected Seifert hypersurface $V$ for $K$. 
   Suppose that $\alpha\in\pi_r(V\x\{t\})$ satisfies the condition $\iota(\alpha)=0$  ($t\in\{1,-1\}$).  

Note that the dimension of $V$ is $n+1$, 
and that $2r\leqq n+1$.  
Hence $\alpha$ is represented by an embedded 
$r$-sphere in $V$. 

Let $\alpha$ also denote this $r$-sphere. 
Then there is a continuous map 
$f:D^{r+1}\to \overline{Q-N(K)}$ 
such that $f(\partial D^{r+1})=\alpha$. 
Since  $\iota(\alpha)=0$, 
 $f($Int$D^{r+1})\cap V=\phi$. 

Note that the dimension of $\overline{Q-N(K)}$ is $n+2$, 
that the dimension of $D^{r+1}$ is $r+1$, 
that $r+1\geqq3$,  
and that $2(r+1)\leqq n+2$. 
Hence we can suppose that $f$ is an embedding map.

Take an embedded $(n+2)$-dimensional $(r+1)$-handle 
whose core is $f(D^{r+1})$ and whose attaching part is embedded in $V$. 
Carry out a surgery on $V$ by using this handle. 
(Note that we remove the interior of $S^r\x D^{n+1-r}$ from $V$ and attach 
$D^{r+1}\x S^{n-r}$.)
We obtain a new $V$. 
Repeating this procedure. 
Since $\pi_r V$ is finitely generated, \newline 
$\iota:\pi_r(V\x\{t\})\to\pi_r(\overline{\overline{Q-N(K)}-N(V)})$ 
becomes injective for $t=1,-1$ after finite times of this procedure. 
This completes the proof of Claim \ref{pi}. \qed

\begin{cla}\label{piinj}
Let $r\leqq[\frac{n-1}{2}]$. 
Let $V$ be an $(r-1)$-connected Seifert hypersurface for $K$.   
Suppose that 
$\pi_r(V\x\{t\})\to\pi_r(\overline{\overline{Q-N(K)}-N(V)})$ is injective for $t=1,-1$. 
Then $\pi_r V=0$. 
\end{cla}

\noindent{\bf Proof of Claim \ref{piinj}.}
Take any element $\alpha\in\pi_r V$. 
Note that the dimension of $V$ is $n+1$ 
and that $2r\leqq n+1$.  
Hence $\alpha$ is represented by an embedded $r$-sphere in $V$. 

Let $\alpha$ also denote this $r$-sphere. 
Since $\pi_r(\overline{Q-N(K)})=0$, 
there is a continuous map $f:D^{r+1}\to\overline{Q-N(K)}$ 
such that $f(\partial D^{r+1})=\alpha$. 

Note that the dimension of $\overline{Q-N(K)}$ is $n+2$, 
that the dimension of $D^{r+1}$ is $r+1$, 
that $r+1\geqq3$,  
and that $2(r+1)\leqq n+2$. 
Hence we can suppose that $f$ is an embedding map. 

We prove the following: 

\begin{cla}\label{sp}
Let $N(\alpha)$ be the tubular neighborhood of $\alpha$ in  $\overline{Q-N(K)}$.   
Then \newline
 $f(D^{p+1})\cap V\cap N(\alpha)=\alpha$.  
\end{cla}

\noindent{\bf Proof of Claim \ref{sp}.}
We need the following claim:

\begin{cla}\label{soku}
Let $p\in\N$. 
Let $q\geqq p+2$. 
Let $\alpha$ be an $\R^q$-bundle over $S^p$. 
Let $\tau$ be the tangent bundle of $S^p$. 
Let $\varepsilon^r$ be the trivial $\R^r$-bundle over $S^p$. 
If $\alpha\oplus\tau=\varepsilon^{p+q}$, then $\alpha=\varepsilon^q$. 
\end{cla}

\noindent{\bf Proof of Claim \ref{soku}.}   
Since $q\geqq p$, we have $\alpha=\varepsilon^1\oplus\beta$, where 
$\beta$ is an $\R^{q-1}$-bundle. 
Since $\tau\oplus\varepsilon^1=\varepsilon^{p+1}$, 
$\alpha\oplus\tau=(\beta\oplus\varepsilon^1)\oplus\tau=
(\tau\oplus\varepsilon^1)\oplus\beta=
\varepsilon^{p+1}\oplus\beta$. 
Hence $\varepsilon^{p+q}
=\varepsilon^{p+1}\oplus\beta$. 

Recall that $\pi_i SO(n)\cong\pi_i(SO(n+1))$ if $1\leqq i\leqq n-2$ and $n\in\N-\{1\}$: 
Reason; The exact sequence $\pi_i SO(n)\to\pi_i(SO(n+1))\to\pi_i S^n.$

Recall that $\R^r$-bundles over $S^p$ are classified by $\pi_{p-1}SO(r)$. 
Note that $p-1\leqq (q-1)-2$ 
and that   $\varepsilon^{p+1}\oplus\beta$ is the trivial $\R^{p+q}$-bundle. 
Hence $\beta$ is the trivial $\R^{q-1}$-bundle over $S^p$. 

Hence $\beta\oplus\varepsilon^1$ is the trivial $\R^q$-bundle over $S^p$. 
This completes the proof of Claim \ref{soku}.\qed\bigbreak

Let $N'(\alpha)$ be the tubular neighborhood of $\alpha$ in $V$ 
(Recall that 
$N(\alpha)$ is the tubular neighborhood of $\alpha$ in  $\overline{Q-N(K)}.$)  
By Claim \ref{soku}, $N'(\alpha)=S^r\x D^{n+1-r}$ and \newline 
$N(\alpha)=S^r\x D^{n+2-r}$. 
Hence $\partial N(\alpha)=S^r\x S^{n+1-r}$ is 
the trivial $S^{n+1-r}$-bundle over $S^r$. 
Since $n+2-r\geqq r+2$,  
all sections of this trivial $S^{n+1-r}$-bundle over $S^r$ are homotopic. 
There is a section that is perpendicular to $N'(V)$ at $\alpha$. 
Another section is defined by $f(D^{r+1})\cap N(\alpha)$. 
Both sections are homotopic.
This completes the proof of Claim \ref{sp}. \qed\bigbreak

Recall that $f$ is an embedding map. 
By Claim \ref{sp} we can suppose that $f(D^{r+1})$ intersects $V$ transversely and that  $f(D^{r+1})\cap V$ is a disjoint union of connected, closed, oriented, $r$-manifolds. 
Note that each of these $r$-manifolds is not an $r$-sphere in general. 
Take an innermost connected $r$-manifold $M$ of these $r$-manifolds. 
Note that $f^{-1}(M)$ in $D^{r+1}$ is diffeomorphic to $M$. 
There is an $(r+1)$-dimensional compact connected, oriented, manifold $W$ 
embedded in $D^{r+1}$ such that $M=\partial W$. 
By the existence of $W$, 
$M$ is a vanishing $r$-cycle in 
$\overline{\overline{Q-N(K)}-N(V)}.$

By Hurewicz theorem 
$\pi_r V=H_r(V;\Z)$. 
By Hurewicz theorem, Mayor-Vietoris theorem, and van Kampen theorem, 
$\pi_1(\overline{\overline{Q-N(K)}-N(V)})=1$, \newline
$\pi_i(\overline{\overline{Q-N(K)}-N(V)})=0$ for $2\leqq i\leqq r-1$, and \newline 
$\pi_r(\overline{\overline{Q-N(K)}-N(V)})= 
H_r(\overline{\overline{Q-N(K)}-N(V)};\Z)$.  
Hence there is 
an $r$-sphere embedded in $V$   
that is homologous to $M$ and 
the homotopy class $[M]\in \pi_r(\overline{\overline{Q-N(K)}-N(V)})$ is zero. 
Since $\pi_r V\to\pi_r(\overline{\overline{Q-N(K)}-N(V)})$ is injective, 
the homotopy class $[M]\in \pi_r V$ is zero. 
By obstruction theory 
there is a continuous map 
$\overline{f}:W\to V$ such that  \newline  
$\overline{f}|_M=f|_M:M\to V.$ 
Hence we remove $M$ from $f(D^{p+1}\cap V)$ and keep the other connected manifolds than $M$ by using a homotopy.

Repeating this procedure, 
we obtain a new $f$ such that $f(D^{p+1}\cap V)=\phi$. 
Hence $\alpha$ is null-homotopic 
in $\overline{\overline{Q-N(K)}-N(V)}$. 
Since $\pi_r V\to\pi_r( \overline{\overline{Q-N(K)}-N(V)} )$ is injective, 
$\alpha$ is null-homotopic in $V$. 
Hence $\pi_r V=0$. 
This completes the proof of Claim \ref{piinj}. \qed\bigbreak

This completes the proof of Proposition \ref{ab}. \qed

\begin{pr}\label{bc}
`$({\rm ii})\Longrightarrow({\rm iii})$' in Theorem $\ref{handle}$ if $n\geqq5$.   
\end{pr}

\noindent{\bf Proof of Proposition \ref{bc}.} 
We prove the following proposition that is stronger than Proposition \ref{bc}. 

\begin{pr}\label{abu}  
Let $n\geqq5$. 
There is a $[\frac{n-1}{2}]$-connected Seifert hypersurface $V$ for $K$.  
Let $B\cup T=K$, where $B$ $($respectively,  $T)$ may be the empty set. 
Let $B$ and $T$ be  
connected, $([\frac{n}{2}]-1)$-connected, compact. 
Then there is a special handle decomposition of $V$  
whose bottom is $B$ and whose top is $T$.  
\end{pr}  

\noindent{\bf Proof of Proposition \ref{abu}.} 
$V$ satisfies the following condition $(*)$: \newline 
There is a handle decomposition  of $V$ with the following properties; 

\smallbreak\noindent(i)   
The bottom is $B$. The top is $T$. 

\smallbreak\noindent(ii) 
It has only one (respectively,  no) 0-handle if $B=\phi$ (respectively,  if $B\neq\phi$). 

\smallbreak\noindent(iii)   
It has only one (respectively,  no) $(n+1)$-handle if $T=\phi$ (respectively,  if $T\neq\phi$).

\smallbreak\noindent 
Reason: 
Use 1-handles and cancel one or some handles if necessary.  

\smallbreak
Note that $V$ is simply-connected and $(n+1)$-dimensional, that $B$ is connected, 
and that $n+1\geqq6$.  
Hence we have the following  (see e.g. \cite[Lemma 1.21 in \S1.3]{Luck}).

\begin{cla}\label{oneh}  
There is a Seifert hypersurface $V$ for the $n$-submanifold $K$ 
whose handle decomposition satisfies the above $(*)$ and  has no 1-handle. 
\end{cla}

\begin{cla}\label{twoh}
There is a Seifert hypersurface $V$ for the $n$-submanifold $K$ 
whose handle decomposition satisfies Claim \ref{oneh} and that has no 2-handle. 
\end{cla}

\noindent{\bf Proof of Claim \ref{twoh}.}  
Take a handle decomposition of $V$ that satisfies Claim \ref{oneh}.  
Take a sub-handle-decomposition   \newline  
\hskip2cm $T_H=$
$\begin{cases} 
\text{(the only one 0-handle)$\cup$(all 2-handles)}& \text{if $B=\phi$} \\  
\text{$(B\x[0,1])\cup$(all 2-handles)}& \text{if $B\neq\phi$} 
\end{cases}$
\newline of this handle decomposition. 

Note that $V$ is parallelizable, that $B$ is simply-connected, 
that the dimension of $B$ is $n$,  
and that $n\geqq5$. 
Hence we have the following condition $(\sharp)$\newline 
\hskip2cm $T_H=$
$\begin{cases}  
\natural^\nu(S^2\x D^{n-1})  &\text{if $B=\phi$}\\  
 (B\x[0,1]\natural^\nu(S^2\x D^{n-1})) &\text{if $B\neq\phi$, }
\end{cases}$  \newline where $\nu\in\{0\}\cup\N$. 

Note that $V$ is connected, $[\frac{n-1}{2}]$-connected, compact,  
and that $B$ is connected, \newline
$([\frac{n}{2}]-1)$-connected, compact. 
Hence $H_2(V,B;\Z)=\pi_2(V,B)=0$ 
by using Mayor-Vietoris theorem, the homotopy exact sequence of pair, Hurewicz theorem. 

By these facts and the above $(\sharp)$ we can eliminate all 2-handles. 
This completes the proof of Claim \ref{twoh}.  \qed

\bigbreak
\noindent{\bf Note.} 
Suppose that two compact connected manifolds $A^a$ intersect $B^b$ transversely in a simply-connected, connected, compact  manifold $C^{a+b+1}$. 
Whitney trick does not work in general when $a$ or $b$ is 2 even if $a+b\geqq5$. 
Reason: 
Let $a$ (respectively,  $b$) be two.  
Whitney disc may intersect $B$ (respectively,  $A$).

\begin{cla}\label{ph}  
There is a Seifert hypersurface $V$ for the $n$-submanifold $K$ 
whose handle decomposition satisfies Claim \ref{twoh}  
and has no $i$-handle $(1\leqq i\leqq[\frac{n-1}{2}])$. 
\end{cla}  

\noindent{\bf Proof of Claim \ref{ph}.}  
Note that $V$ is connected, $[\frac{n-1}{2}]$-connected, compact,   
and that $B$ is connected, $([\frac{n}{2}]-1)$-connected, compact. 
Hence 
$H_i(V, B;\Z)=\pi_i (V, B)=0  \newline (1\leqq i\leqq[\frac{n-1}{2}])$ 
by using Mayor-Vietoris theorem, the homotopy exact sequence of pair, Hurewicz theorem. 
By this fact and $n+1\geqq6$,  
we can eliminate all $i$-handles ($1\leqq i\leqq[\frac{n-1}{2}]$) 
by using Whitney trick. 
This completes the proof of Claim \ref{ph}. \qed\bigbreak

Take the dual handle decomposition of the handle decomposition that satisfies Claim \ref{ph}. Eliminate one or some handles if necessary 
in the same manner as above.   
This completes the proof of Proposition \ref{abu}. \qed\bigbreak

This completes the proof of Proposition \ref{bc}. \qed\bigbreak


\begin{pr}\label{nisansan}
`$({\rm ii})\Longrightarrow({\rm iii})$'  in Theorem $\ref{handle}$ if $n=3,4$.   
\end{pr}

\noindent{\bf Proof of Proposition \ref{nisansan}.} 
It is trivial that 
Proposition \ref{nisansan} follows from the following proposition. 

\begin{pr}\label{abu34}
The $n=3,4$ case of Proposition $\ref{abu}$ holds. 
\end{pr}

\noindent{\bf Proof of Proposition \ref{abu34}.} 
Take a simply-connected Seifert hypersurface $V$ for $K\newline\subset S^{n+2}$. 
Since $V$ is connected, compact, and $B$ is connected, 
$V$ satisfies the condition $(*)$ in the first paragraph of Proof of Proposition \ref{abu}. 
Take a 1-handle $h^1$ of the handle decomposition. 
Since $V$ is oriented, $B\cup h^1=B\natural(S^1\x D^n)$. 
Since $\pi_1V=1$, 
there is a continuous map $f:D^2\to V$ such that 
$f(\partial D^2)=S^1\x\{0\}$. 
Push off $f({\rm Int} D^2)$ from $V$, keeping $f(\partial D^2)$,  
into the positive direction of the normal bundle of 
$V$ in $\overline{S^{n+2}-N(K)}$.   
Thus we obtain a continuous map $g:D^2\to S^{n+2}-N(K)$ 
such that $D^2\cap V=\partial D^2=S^1\x\{0\}$.

Note that the dimension of $\overline{S^{n+2}-N(K)}$ is $n+2$, 
that the dimension of $D^2$ is 2, 
and that $n+2\geqq 2\x2$. 
Hence we can suppose that $g$ is an embedding map.

Take an embedded $(n+2)$-dimensional 2-handle whose core is $D^2$ and 
whose attaching part is embedded in $V$.  
Carry out a surgery on $V$ by using this 2-handle. 
(Note $n-1\neq1$. Note that we remove the interior of $S^1\x D^n$ from $V$ and add $D^2\x S^{n-1}$.)
Then the 1-handle $h^1$ is eliminated and a new $(n-1)$-handle is obtained. 
(Note that the dual of an $(n-1)$-handle is a 2-handle not a 1-handle.) 
We eliminate the 1-handle $h^1$ and obtain a new $V$. 
Repeating this procedure, we eliminate all 1-handles.

Take its dual handle decomposition.   
Since $V$ is connected, compact, and $B$ is connected, 
it also satisfies the condition $(*)$ 
(Of course we replace $B$ (respectively,  $T$) in the condition $(*)$ with  $T$ (respectively,  $B$)).  
It has no $n$-handle.   
We can eliminate all 1-handles as above so that we do not obtain a new $n$-handle. 

Therefore the new $V$ has a handle decomposition 
such that the bottom is $B$, that the top is $T$, 
that it satisfies the condition $(*)$, 
and that 
it has no 1-handle, no $n$-handle.
This completes the proof of Proposition \ref{abu34}. \qed\bigbreak

This completes the proof of Proposition \ref{nisansan}. \qed\bigbreak

This completes the proof of Theorem \ref{handle}. \qed\bigbreak

\section{Proof of Theorem \ref{Tokyo}}\label{sdochi} 
\noindent
We prove Theorem \ref{Tokyo}.(ii)$\Longrightarrow$Theorem \ref{Tokyo}(i): 
Attach an embedded $(2p+2)$-dimensional $(p+1)$-handle 
to a simple Seifert hypersurface for $K$,  
carry out a surgery by using this handle, and obtain a new one.
Then a simple 
Seifert matrix associated with the old one is $(-1)^p$-S-equivalent to 
that associated with the new one.  
Repeat this procedure. 

\begin{cla}\label{lakehigh}
 Theorem $\ref{Tokyo}.{(\rm i)}$  $\Longrightarrow$ Theorem $\ref{Tokyo}.{(\rm ii)}$.   
\end{cla}

\noindent{\bf Proof of Claim \ref{lakehigh}.} 
Let $(X,Y)$ denote a pair (A manifold, its submanifold).   
Let $(X,Y)\x[0,1]$ denote a pair 
(The manifold$\x[0,1]$, its submanifold$\x[0,1]$) which is a level preserving embedding.   
Take $(S^{2p+3}, K)\x[0,1]$. 
Let $V$ (respectively,  $V'$) 
be a simple Seifert hypersurface for $K$ 
whose simple 
Seifert matrix is $P$ (respectively,  $P'$). 
Take a $(2p+2)$-submanifold $V\cup (K\x[0,1])\cup V'$ 
in a $(2p+4)$-manifold $S^{2p+3}\x[0,1]$. 
By Propositions \ref{abu} and \ref{abu34} 
we can suppose that 
a $(2p+3)$-dimensional Seifert hypersurface $W$ for 
the $(2p+2)$-submanifold $V\cup (K\x[0,1])\cup V'$ 
has a special handle decomposition \newline
$(V\x[0,1])\cup$($(p+1)$-handles)$\cup$($(p+2)$-handles)$\cup(V'\x[0,1])$. 
Let $\Phi:W\to[0,1]$ be a height function and a Morse function  
which gives this handle decomposition. 
By using isotopy we can suppose the following: 
Let $0=t_0\leqq t_1...\leqq t_\nu=1$ be a partition of $I$ satisfying 

\smallbreak\noindent(i) 
Each $t_i$ is a regular value of $\Phi$. 

\smallbreak\noindent(ii) 
At most one critical value of $\Phi$ lies in each interval $(t_i, t_{i+1})$. 

\smallbreak 
Note that for each $i$, $\Phi^{-1}(t_i)$ is connected and $p$-connected. 

\smallbreak 
Therefore it suffices to prove the case where $\Phi$ has only one critical point. If we change $[0,1]$ to $[0,-1]$, then the index $\xi$ of critical point 
becomes $2p+3-\xi$. Hence   
it suffices to prove the case where $\Phi$ has only one critical point of index $(p+1)$. 
 That is, it suffices to prove the case where  
$W$ has a handle decomposition \newline
$(V\x[0,1])\cup$(only one $(p+1)$-handle)$\cup(V'\x[0,1])$. 
We can suppose that 
this only one $(p+1)$-handle is attached to $V$ in $S^{2p+3}\x\{t\}$ for a $t$ by a famous method 
in `Morse theory with handle bodies'.  
Hence a simple 
Seifert matrix for $V$ is $(-1)^p$-$S$-equivalent to that for $V'$. 
This completes the proof of Claim \ref{lakehigh}. \qed\bigbreak

This completes the proof of Theorem \ref{Tokyo}.\qed\bigbreak

\section{Proof of Proposition \ref{Osaka}} \label{ichio}
\noindent
We show an example. 
Embed $S^3\x D^3$ in $S^7$ such that $S^3\x$$\{$the center$\}$ is 
the standard 3-sphere trivially embedded in $S^7$.  
Let the submanifold $\partial(S^3\x D^3)$ in $S^7$ satisfy 
that 
the $S^3\x D^3$ is its simple Seifert hypersurface   
whose simple 
Seifert matrix is a $1\x1$-matrix $(0)$.  
This 5-submanifold is called $K$.  
We can suppose that $D^4\x S^2$ bounds $K$.  
Hence $\phi$ and 
$\begin{pmatrix}
0&1\\
0&0
\end{pmatrix}$ are 3-Seifert matrices for $K$, 
where $\phi$ is the empty matrix.  
Note that 
$(0)$ is not $S$-equivalent to 
$\phi$ or 
$\begin{pmatrix}
0&1\\
0&0
\end{pmatrix}$.

We show another example. 
There is a closed oriented 3-submanifold $K$ of $S^5$ with the following properties.  
(See \cite[\S10]{KauffmanOgasa} for this submanifold for detail.)  

\smallbreak\noindent
(1)
There is a Seifert hypersurface $V$ for $K$ whose framed link representation is 
the $(2, 2a)$ torus link such that  the framing of each component is zero 
($a\in\N-\{1\}$). 
$\begin{pmatrix}
0&a\\
0&0
\end{pmatrix}$ is a simple 
Seifert matrix for $K$. 

\smallbreak\noindent
(2) 
There is a Seifert hypersurface $W$ for $K$ whose framed link representation is 
the $(2, 2a)$ torus link such that  the framing of one component is zero and the other component is the dot circle. (Carry out a surgery on $V$ by using a 5-dimensional 3-handle 
embedded in $S^5$ and obtain $W$.) $W$ is a 4-dimensional homology ball. 
Hence $\phi$ and 
$\begin{pmatrix}
0&1\\
0&0
\end{pmatrix}$ are Seifert matrices for $K$.  
Note that 
$\begin{pmatrix}
0&a\\
0&0
\end{pmatrix}$, where $a\in\N-\{1\}$, 
is not $(-1)$-$S$-equivalent to 
$\phi$ or 
$\begin{pmatrix}
0&1\\
0&0
\end{pmatrix}$. 
\qed

\bigbreak
\section{Proof of Theorem \ref{Chicago}} \label{Seiequ}
\noindent
First we prove Theorem \ref{Chicago}.(1). 
Theorem \ref{Tokyo} implies that `$({\rm i})\Longrightarrow({\rm ii})$' in Theorem \ref{Chicago}.$(1)$.

\begin{lem}\label{river}
`$({\rm ii})\Longrightarrow({\rm i})$' in Theorem $\ref{Chicago}.(1)$.  
\end{lem}
 
\noindent{\bf Proof of Proposition \ref{river}.} 
It suffices to prove the case 
where $P_J=P_K$.  
By Theorem \ref{Tokyo} and `Theorem \ref{handle} and its proof',  we can suppose that 
there is a simple Seifert hypersurface  $V_*$ 
with a special handle decomposition (the top $\partial V_*$, the bottom $\phi$) 
whose simple 
Seifert matrix is $P_*$ ($*=J,K$). 
By Proposition \ref{daiji}, 
$P_*+(-1)^p\hskip1mm^t\hskip-1mm P_* (*=J, K)$ is the intersection product on 
$H_p(V_*;\Z)$. 

Embed $(2p+2)$-dimensional 0-handle $h^0$ in $S^{2p+3}$ trivially.  
We can take $p$-spheres in $\partial h^0$ 
and attach embedded $(2p+2)$-dimensional $(p+1)$-handles 
to $h^0$ along the $p$-spheres so that a simple Seifert matrix associated with 
the result $h^0\cup$($(p+1)$-handles) is $P_*$ ($*=J,K$). 

\bigbreak

Then $V_J$ is diffeomorphic to $V_K$. Reason: 
The core of the attached part of each $(2p+2)$-dimensional $(p+1)$-handle 
is a $p$-sphere. The boundary of $h^0$ is a $(2p+1)$-sphere. 
Furthermore $p\geqq2$. 
By \cite{Haefligerunknot, HaefligerFr, Whitney, Whitneytrick, Wu}, 
the embedding type of an ordered   disjoint  union of $p$-spheres in $\partial h^0$ 
is determined by the set of the linking numbers of each pair of $p$-spheres. 

\bigbreak
The core of each $(2p+2)$-dimensional $(p+1)$-handle is $(p+1)$-dimensional. $S^{2p+3}$ is $(2p+3)$-dimensional. 
Furthermore $p\geqq2$. 
By \cite{Haefligerunknot, HaefligerFr, Whitney, Whitneytrick, Wu},  
the embedding type of an ordered  disjoint union of $(p+1)$-handles in $S^{2p+3}-$Int$h^0$ keeping the attached part of each $(p+1)$-handle 
is determined by 
the Seifert matrix. 

Therefore there is a diffeomorphism map $f:S^{2p+3}\to S^{2p+3}$ such that 
$f(V_J)=V_K$. 

This completes the proof of Proposition \ref{river}.\qed

\bigbreak 
Next we prove Theorem \ref{Chicago}.(2). 
`$({\rm i})\Longrightarrow({\rm ii})$' in Theorem \ref{Chicago}.$(2)$ is trivial.

\begin{lem}\label{kawa}
`$({\rm ii})\Longrightarrow({\rm i})$' in Theorem $\ref{Chicago}.(2)$.  
\end{lem}

\noindent{\bf Proof of Proposition \ref{kawa}.} 
By \cite{Haefligerunknot, HaefligerFr, Whitney, Whitneytrick, Wu}, 
the embedding type of an ordered  disjoint union of 2-handles in $S^5-$Int$h^0$ keeping the attached part of each 3-handle 
is determined by 
the Seifert matrix.

Therefore submanifolds, $V_J$ and $V_K$, of $S^5$ are isotopic.   
This completes the proof of Proposition \ref{kawa}.\qed\bigbreak

This completes the proof of Theorem \ref{Chicago}. \qed

\bigbreak
\section{Proof of Theorem \ref{aoiro}}\label{blue}
\noindent 
In this section and the following one 
we use relations between 
Seifert matrices and knot products. 
See \cite[Proposition 5.4]{KauffmanOgasa} for the relations. 
See \S\ref{RkB} of this paper for $S_*$ and $N_\#$ in the proposition. 
See also \cite[section 6]{KauffmanNeumann}, 
and 
\cite[page 541]{LevineAlex}.

\begin{cla}\label{ringo}
There is a simple 
Seifert matrix $P_*$ for $*$ $(*=J,K)$ with the following property:  Each element of $P_J$ is the same as that of $P_K$ 
except for only one diagonal element.  They differ by one. 
\end{cla}

\noindent{\bf Proof of Claim \ref{ringo}.}
There is a $(2p+3)$-ball $B$ where the twist-move is carried out. 
Take a $(2p+2)$-dimensional $(p+1)$-handle $h$ associated with the twist move that is embedded in $B$.    
Let $Z$ be a $(2p+1)$-dimensional closed oriented submanifold 
 $(J-B)\cup (h\cap \partial B)\subset S^{2p+3}$.  
Thus $(J,K,Z)$ is a twist-move triple. 
Note that $Z$ is embedded in $\overline{S^{2p+3}-B}$.  
By the construction of $Z$, 
$Z$ is $(p-1)$-connected. Furthermore 
 $\overline{S^{2p+3}-N(J)}$
(respectively,   $\overline{S^{2p+3}-N(K)}$)
is made by attaching one $(2p+3)$-dimensional $(p+2)$-handle
to $\overline{\overline{S^{2p+3}-B}-N(Z)}$. 
Hence 
 $\pi_i(\overline{S^{2p+3}-N(J)})$
 $=\pi_i(\overline{S^{2p+3}-N(K)})$
 $=\pi_i(\overline{\overline{S^{2p+3}-B}-N(Z)})$ 
($1\leqq i\leqq p$). 

By the assumption, $J=A\otimes[2]$ holds and $A$ is 
a $(2p-1)$-dimensional connected, $(p-2)$-connected, simple submanifold of $S^{2p+1}$. 
Hence, by \cite{Kauffman, KauffmanNeumann}, 
$J$ is 
a $(2p+1)$-dimensional connected, $(p-1)$-connected, simple submanifold of $S^{2p+3}$.

Hence $Z$ is 
a $(2p+1)$-dimensional connected, $(p-1)$-connected, simple 
submanifold of $\overline{S^{2p+3}-B}$.   
By Theorem \ref{handle} there is 
a $p$-connected Seifert hypersurface $W$ in $\overline{S^{2p+3}-B}$ for $Z$ 
that has a special handle decomposition (the top $Z$, the bottom $\phi$). 
Note that the handle $h$ is attached to $W$. 
By using this $W\cup h$ we obtain 
a simple Seifert hypersurface $V_J$ (respectively,  $V_K$) for $J$ (respectively,  $K$)
that has a special handle decomposition. 
This completes the proof of Claim \ref{ringo}. \qed

\begin{cla}\label{fuyu}
$(-1)^{p-1}P_J$ is a simple 
Seifert matrix for $A$. 
\end{cla}

\noindent{\bf Proof of Claim \ref{fuyu}.}
Let $P_A$ be a simple 
Seifert matrix for $A$.$\cdot\cdot\cdot(*)$   \newline
We prove the following Claims \ref{mikan} and \ref{pear}. 

\begin{cla}\label{mikan}
$(-1)^{p-1}P_A$ is $(-1)^p$-$S$-equivalent to $P_J$. 
\end{cla}

\noindent{\bf Proof of Claim \ref{mikan}.}
Since $A\otimes[2]=J$, 
$(-1)^{p-1}P_A$ is a simple 
Seifert matrix for $J$. 
(Reason: Note \ref{MA}. 
Note that 
a Seifert matrix for $[2]$ is 
a $1\x1$-matrix $\Lambda_2=(1)$.  
Recall that $\Lambda_2$ is defined right before Theorem \ref{KNB}.)
By Theorem \ref{Chicago}, we have Claim \ref{mikan}. \qed

\begin{cla}\label{pear}
 $(-1)^{p-1}P_A$ is $(-1)^{p+1}$-$S$-equivalent to $P_J$.
\end{cla}

\noindent{\bf Proof of Claim \ref{pear}.}
By Note \ref{MA},  
$(-1)^p P_J$ is a simple 
Seifert matrix for $J\otimes[2]$.  
By Note \ref{MA} and Claim \ref{mikan},  
$(-1)^p(-1)^{p-1}P_A$ is a simple 
Seifert matrix for $J\otimes[2]$.  
By Theorem \ref{Chicago}  
 it holds that  $(-1)^{p}(-1)^{p-1}P_A$ is $(-1)^{p+1}$-$S$-equivalent to 
$(-1)^{p}P_J$.
Hence we have Claim \ref{pear}. \qed

\bigbreak\noindent{\bf Note.}
Claims \ref{mikan} and \ref{pear} hold on time. 
 Reason: Recall Proposition \ref{daiji}. 
The way of making the intersection matrix 
from a pair of related positive Seifert matrix and negative one depends on 
whether $p$ is odd or even in general. 
\bigbreak

By Claim \ref{pear}  $P_A$ is $(-1)^{p+1}$-$S$-equivalent to $(-1)^{p-1}P_J$.
This fact, the above $(*)$, and Theorem \ref{Tokyo} imply 
Claim \ref{fuyu}.   \qed\bigbreak

By Claims \ref{ringo} and \ref{fuyu}, 
we can make a $(2p-1)$-dimensional connected, $(p-1)$-connected, simple 
submanifold $B$ of $S^{2p+1}$ from $A$ by one twist move 
so that a simple 
 Seifert matrix of $B$ is $(-1)^{p-1}P_K$. 
By \cite{Kauffman, KauffmanNeumann} and Theorem \ref{Chicago}, 
we have Theorems \ref{aoiro}.(i), (ii), and (iii). 
\qed

\bigbreak
\section{Proof of Theorem \ref{Komagome}}\label{honkomagome}
\begin{cla}\label{umi}
There is a simple 
Seifert matrix $P_*$ for $*$ $(*=J,K)$ with the following property: 
Let $p_{*, ij}$ be the $(i,j)$-element of $P_*$. 
Let $P_*$ be a $c\x c$-matrix $(c\in\N)$. 
There are natural numbers $a, b$  $\leqq c$ such that  $a\neq b$ 
and that 
$$
\left\{
\begin{array}{ll}
p_{J,ij}=p_{K,ij}-1 & \text{ if $(i,j)=(a, b)$}\\
p_{J,ij}=p_{K,ij}    & \text{ if $(i,j)\neq(a, b)$ and if $(i,j)\neq(b, a). $}\\
\end{array}
\right.
$$
\noindent 
$($Recall that the $(b, a)$-element is the same as the $(a, b)$-element.$)$   
\end{cla}

\noindent{\bf Proof of Claim \ref{umi}.}
There is a $(4\mu+3)$-ball $B$ where the $(2\mu+1, 2\mu+1)$-pass-move is carried out. 
Take $(4\mu+2)$-dimensional $(2\mu+1)$-handles $h$ and $h'$ 
associated with the $(2\mu+1, 2\mu+1)$-pass-move that are embedded in $B$.    
Let $Z$ be a $(4\mu+1)$-dimensional closed oriented submanifold 
 $(J-B)\cup(h\cap\partial B)\cup(h'\cap\partial B) \subset S^{4\mu+3}$.  
Thus $(J,K,Z)$ is a \newline 
$(2\mu+1, 2\mu+1)$-pass-move triple. 
Note that $Z$ is embedded in $\overline{S^{4\mu+3}-B}$.  
By the construction of $Z$,   
$Z$ is $(2\mu-1)$-connected. Furthermore  
 $\overline{S^{4\mu+3}-N(J)}$
(respectively,   $\overline{S^{4\mu+3}-N(K)}$)
is made by attaching 
 two $(4\mu+3)$-dimensional $(2\mu+2)$-handles 
and one $(4\mu+3)$-dimensional $(4\mu+2)$-handles   
to 
$\overline{\overline{S^{4\mu+3}-B}-N(Z)}$.
Hence 
 $\pi_i(\overline{S^{4\mu+3}-N(J)})$
 $=\pi_i(\overline{S^{4\mu+3}-N(K)})$
 $=\pi_i(\overline{\overline{S^{4\mu+3}-B}-N(Z)})$ 
($1\leqq i\leqq2\mu$). 

By the assumption $J=A\otimes^\mu$Hopf holds and $A$ is a 1-link. 
By \cite{Kauffman, KauffmanNeumann}  
$J$ is 
a $(4\mu+1)$-dimensional connected, $(2\mu-1)$-connected, simple 
submanifold of $S^{4\mu+3}$. 

Hence $Z$ is 
a $(4\mu+1)$-dimensional connected, $(2\mu-1)$-connected, simple submanifold 
of $\overline{S^{4\mu+3}-B}$.   
By Theorem \ref{handle} there is 
a $2\mu$-connected Seifert hypersurface $W$ in $\overline{S^{4\mu+3}-B}$ for $Z$ 
that has a special handle decomposition  (the top $Z$, the bottom $\phi$). 
Note that the handles $h$ and $h'$ are attached to $W$. 
By using this $W\cup h\cup h'$ we obtain 
a simple Seifert hypersurface $V_J$ (respectively,  $V_K$) for $J$ (respectively,  $K$) 
that has a special handle decomposition. 
Hence we have Claim \ref{umi}. \qed\bigbreak

Let $P_A$ be a Seifert matrix for the 1-knot $A$.$\cdot\cdot\cdot(\#)$ 

\begin{cla}\label{shiwasu}
$P_A$ is $S$-equivalent to $(-1)^\mu P_J$. 
\end{cla}

\noindent{\bf Proof of Claim \ref{shiwasu}.}
Since $J=A\otimes^\mu$Hopf, 
 $(-1)^\mu P_A$ is 
a simple 
Seifert matrix for $J$.  
(Reason: See Note \ref{MA}. 
Note that 
a Seifert matrix for the Hopf link is 
a $1\x1$-matrix $\Lambda_{2,2}=(-1)$.   
Recall that $\Lambda_{2,2}$ is defined right before Theorem \ref{KNB}.)
By Theorem \ref{Chicago}   
$(-1)^\mu P_A$ is $S$-equivalent to $P_J$. 
Hence we have Claim \ref{shiwasu}. \qed\bigbreak

By Claim \ref{shiwasu} we have the following: 
\begin{cla}\label{nenmatsu}
There is a 1-link $A'$ whose Seifert matrix is $(-1)^\mu P_J$. 
\end{cla}

By Claim \ref{umi} and  \ref{nenmatsu}, we have the following:  

\begin{cla}\label{nenshi} 
We can make a 1-link $B$ from $A'$ by one pass-move 
so that a Seifert matrix of $B$ is $(-1)^\mu P_K$.
\end{cla}

$A$ is pass-equivalent to $A'$ 
because of the above $(\#)$,  Claims  \ref{shiwasu} and \ref{nenmatsu}.   
By this and Claim \ref{nenshi}  $A$ is pass-equivalent to $B$. 
By Theorem \ref{Chicago} and Claim \ref{nenshi} 
$B\otimes^\mu$Hopf is isotopic to $K$.  
This completes the proof of Theorem \ref{Komagome}.\qed

\begin{note}\label{ei}
\cite[Main Theorem 4.1]{KauffmanOgasaII}, which is cited right above Theorem \ref{Bunkyo}, 
is also proved by using Theorem \ref{Chicago}. Outline of the proof: 
There is a Seifert surface 
$V_A$ (respectively,  $V_B$) 
for   
$A$ (respectively,  $B$) 
such that $V_A$ is obtained from $V_B$ by one pass-move. 
There is a Seifert matrix $S_A$ (respectively,  $S_B$) associated with $V_A$ (respectively,  $V_B$). 
By \cite{KauffmanNeumann},  
$(-1)^\mu S_A$ (respectively,  $(-1)^\mu S_B$) is a $(2\mu+1)$-Seifert matrix for 
$K_A\otimes^\mu$Hopf (respectively,  $K_B\otimes^\mu$Hopf).  
By the definition of the $(2\mu+1, 2\mu+1)$-pass-move and Theorem \ref{Chicago}, 
there is a simple Seifert hypersurface 
$W_A$ (respectively,  $W_B$) for $A\otimes^\mu$Hopf (respectively,  $B\otimes^\mu$Hopf) 
such that 
$(-1)^\mu S_A$ (respectively,  $(-1)^\mu S_B$) is a $(2\mu+1)$-Seifert matrix 
associated with $W_A$ (respectively,  $W_B$) 
and 
such that $W_A$ is obtained from $W_B$ by one  $(2\mu+1, 2\mu+1)$-pass-move.  
Hence \cite[Main Theorem 4.1]{KauffmanOgasaII} holds. 
\end{note}

\noindent
{\footnotesize 
{\bf Acknowledgment.}  
The authors would like to thank Osamu Saeki for the valuable discussion.  
}

\noindent
Louis H. Kauffman: 
Department of Mathematics, Statistics, and Computer Science,  
University of Illinois at Chicago, 
851 South Morgan Street, 
Chicago, Illinois 60607-7045, USA  \quad
kauffman@uic.edu
\\

\noindent
Eiji Ogasa:  Computer Science, Meijigakuin University, Yokohama, Kanagawa, 244-8539, Japan 
\quad pqr100pqr100@yahoo.co.jp  \quad
ogasa@mail1.meijigkakuin.ac.jp

\end{document}